\documentclass[11pt]{article}

\usepackage[slantedGreek]{mathpazo}
\usepackage[pdftex]{graphicx}
\usepackage{amsfonts}
\usepackage{setspace}
\usepackage{amssymb}
\usepackage{amsthm}
\usepackage{graphicx}
\usepackage{amsmath}
\usepackage{mathtools}
\usepackage{lscape}
\usepackage{subcaption}
\usepackage{suffix}
\usepackage{color}

\usepackage[sort,comma]{natbib}

\allowdisplaybreaks

\usepackage{suffix}

\newcommand{\E}{{\bf E}}

\newcommand{\Nor}{{\cal N}}  

\renewcommand{\P}{{\bf P}}

\newcommand{\V}{{\bf V}}

\newcommand{\ben}{\begin{enumerate}}
\newcommand{\een}{\end{enumerate}}
\newcommand{\beq}{\begin{equation}}
\newcommand{\eeq}{\end{equation}}

\newcommand{\half}{\frac{1}{2}}

\newcommand{\norm}[1]{\lVert#1\rVert}

\newcommand{\ba}{{\bf a}}
\newcommand{\bb}{{\bf b}}
\newcommand{\bc}{{\bf c}}
\newcommand{\bd}{{\bf d}}

\newcommand{\defeq}{\operatorname{:=}}

\newcommand{\estHs}{\ensuremath{{\hat{\theta}}_{HS}}{}}
\newcommand{\estJs}{\ensuremath{{\hat{\theta}}_{JS}}{}}
\newcommand{\MSE}{\mathrm{MSE}}
\DeclarePairedDelimiterX\MeijerM[3]{\lparen}{\rparen}%
{\begin{smallmatrix}#1 \\ #2\end{smallmatrix}\delimsize\vert\,#3}

\newcommand\MeijerG[8][]{%
  G^{\,#2,#3}_{#4,#5}\MeijerM[#1]{#6}{#7}{#8}}

\WithSuffix\newcommand\MeijerG*[7]{G^{\,#1,#2}_{#3,#4}\MeijerM*{#5}{#6}{#7}}

\DeclarePairedDelimiterX\pFqM[3]{\lparen}{\rparen}%
{\begin{smallmatrix}#1 \\ #2\end{smallmatrix}\delimsize\vert\,#3}

\newcommand\pFq[6][]{%
  {}_{#2}F_{#3}\pFqM[#1]{#4}{#5}{#6}}

\WithSuffix\newcommand\pFq*[5]{{}_{#1}F_{#2}\pFqM*{#3}{#4}{#5}}

\newtheorem{theorem}{THEOREM}
\numberwithin{theorem}{section}

\newtheorem{Def}{DEFINITION}
\numberwithin{Def}{section}
\newtheorem{remark}{REMARK}
\numberwithin{remark}{section}

\newtheorem{proposition}{PROPOSITION}
\numberwithin{proposition}{section}

\numberwithin{lemma}{section}

\numberwithin{Cor}{section}

\newtheorem{assn}{ASSUMPTION}

\usepackage{booktabs,array}

\newcount\rowc

\makeatletter
\def\ttabular{%
\hbox\bgroup
\let\\\cr
\def\rulea{\ifnum\rowc=\@ne \hrule height 1.0pt \fi}
\def\ruleb{
\ifnum\rowc=1\hrule height 1.0pt  \else
\ifnum\rowc= 3  \hrule height 0.5pt \else
\ifnum\rowc= 5  \hrule height 0.5pt \else
\ifnum\rowc= 7  \hrule height 0.5pt \else
\ifnum\rowc= 9  \hrule height 0.5pt \else
\ifnum\rowc= 11  \hrule height 0.5pt 
  \else \hrule height 0pt
\fi\fi\fi\fi\fi\fi}
\valign\bgroup
\global\rowc\@ne
\rulea
\hbox to 7em{\strut \hfill##\hfill}%
\ruleb
&&%
\global\advance\rowc\@ne
\hbox to 7em{\strut\hfill##\hfill}%
\ruleb
\cr}
\def\endttabular{%
\crcr\egroup\egroup}

\graphicspath{{../../figures/}{../figures/}{./figures/}{./}}
\renewcommand{\E}{\mathbb E}
\renewcommand{\P}{\mathbb P}
\renewcommand{\V}{\mathbb V}

\numberwithin{equation}{section}

\title{The Horseshoe+ Estimator of Ultra-Sparse Signals}
\date{}
\begin{document}
\maketitle
\baselineskip=15pt
\begin{center}
\vspace{-1cm}
Anindya Bhadra\\
Department of Statistics, Purdue University, 250 N. University Street, West Lafayette, IN 47907-2066\\
bhadra@purdue.edu\\
\hskip 5mm \\
Jyotishka Datta\\
Statistical and Applied Mathematical Sciences Institute, 19 T.W. Alexander Drive,  Research Triangle Park, NC 27709-4006.\\
Department of Statistical Science, Duke University \\
jd298@stat.duke.edu\\
\hskip 5mm \\
Nicholas G. Polson and Brandon Willard\\
The University of Chicago Booth School of Business, 5807 S. Woodlawn Ave., Chicago, IL 60637\\
ngp@chicagobooth.edu, brandonwillard@gmail.com\\
\end{center}

%
%
%
%
%
%
%

\begin{abstract}

\noindent We propose a new prior for ultra-sparse signal detection that we term the ``horseshoe+ prior.'' The horseshoe+
prior is a natural extension of the horseshoe prior that has
achieved success in the estimation and detection of sparse signals and has been shown to possess a number of desirable theoretical properties while enjoying computational feasibility in high dimensions.   The horseshoe+ prior builds upon these advantages. Our work 
proves that the horseshoe+ posterior concentrates at a rate faster than that of the horseshoe in the Kullback-Leibler (K-L) sense. We also establish theoretically that the proposed estimator has lower posterior mean squared error in estimating signals compared to the horseshoe and achieves the optimal Bayes risk in testing up to a constant. For global-local scale mixture priors, we develop a new technique for analyzing the marginal sparse prior densities using the class of Meijer-G functions. In simulations, the horseshoe+ estimator demonstrates superior performance in a standard design setting against competing methods, including the horseshoe and Dirichlet-Laplace estimators. We conclude with an illustration on a prostate cancer data set and by pointing out some directions for future research.\\

\end{abstract}
%



\section{Introduction}

Ultra-sparse signal detection provides a challenge for developing statistical estimators. In the classical normal means inference problem, we observe data from the probability model $  (y_i | \theta_i)  \sim \Nor ( \theta_i,1)$ for $i = 1, \ldots, n$. We wish to provide an estimator for the vector of  normal means $ \theta = ( \theta_1, \ldots , \theta_n ) $. 
Sparsity occurs when a large portion of the parameter vector contains zeros.  The ``ultra-sparse'' or ``nearly black'' vector case occurs when the parameter vector $\theta$ lies in the set $ l_0 [ p_n] \equiv \{ \theta : \# ( \theta_i \neq 0 ) \leq p_n \} $ with the number of non-zero parameter values $ p_n = o(n) $ where $ p_n \rightarrow \infty$ as $ n \rightarrow \infty $.  

To motivate the need for developing new prior distributions, consider 
the classic James-Stein ``global'' shrinkage rule, $\estJs(y)$. This estimator uniformly dominates the traditional sample mean estimator, $\hat{\theta}$. For all values of the true parameter $\theta$ and for $n>2$, we have the classical mean squared error (MSE) risk bound:
$$ 
R(\estJs, \theta) \defeq \E_{y|\theta} {\Vert \estJs(y) - \theta \Vert}^2 < n 
    = \E_{y|\theta} {\Vert y - \theta \Vert}^2, \;\;\; \forall \theta.
$$
However, for a sparse signal, $\estJs(y)$ performs poorly. Suppose that the true parameter $\theta$ is an
``$r$-spike'' with $r$ coordinates of magnitude $ \sqrt{n/r}$ and the rest set at zero, giving $ {\Vert \theta \Vert}^2 =n$.
Then \cite{johnstone2004needles} showed that the classical risk satisfies 
$ R \left ( \estJs , \theta \right ) \geq n/2 $ whereas simple thresholding at 
$ \sqrt{2 \log n}$ performs with risk $\sqrt{\log n }$. 

To address this issue, a ``global-local'' shrinkage estimator called the horseshoe estimator was proposed by \cite{carvalho2010horseshoe}. The horseshoe estimator, $\estHs(y)$, provides a Bayes rule that inherits good MSE properties of global shrinkage estimators and simultaneously provides asymptotic minimax risk for estimating sparse signals.  For example, \cite{polson2012half} showed that $\estHs(y)$ uniformly dominates the traditional sample mean estimator in terms of MSE and \cite{van2014horseshoe} showed that the horseshoe estimator has good posterior concentration properties.
Specifically, the horseshoe estimator achieves
$$ 
\sup_{ \theta \in l_0[p_n] } \; 
\mathbb{E}_{ y | \theta } \norm{ \estHs (y) - \theta }^2 \asymp
p_n \log \left ( n / p_n \right ),
$$
which is the asymptotically minimax risk rate in $\ell_2$ for nearly black objects \citep{donoho1992maximum}.
Here the ``worst'' $\theta \in l_0[p_n]$ is obtained at the maximum absolute difference
$\left| \estHs(y) - y \right|$ where $\estHs(y) = \mathbb{E}_{HS}(\theta|y)$ can be
interpreted as a Bayes posterior mean which is optimal under the Bayes MSE. 

Though the horseshoe prior was originally designed to provide an accurate and efficient estimator of a sparse normal mean vector, it turns out that the multiple testing rule induced by the horseshoe prior also enjoys the ``oracle property'' in testing under the 0-1 loss \citep{datta2013asymptotic}. For the multiple testing problem in the classical two-groups model, many approaches involve explicitly modeling the ultra-sparse mean as a mixture of a point mass at zero and a heavy-tailed alternative, also known as the ``spike-and-slab'' approach \citep{mitchell88}. This results in a posterior distribution over a high-dimensional discrete space, exploring which often leads to extreme computational cost. The one-group model, inspired by the widespread popularity of the lasso for variable selection in regression \citep{tibs96}, is computationally more tractable, and can be used to select a model through concentration of measure in a space of pseudo-probabilities, rather than in the $n$-dimensional Euclidean space \citep{carvalho2010horseshoe,polson2010shrink,datta2013asymptotic}. In particular, the horseshoe prior leads to ``pseudo-posterior" probabilities that mimic the true posterior inclusion probabilities from a two-groups mixture model, and induces a multiple testing rule with attractive properties. 
Specifically, \cite{datta2013asymptotic} proved that the Bayes risk for the horseshoe estimator attains the Bayes risk of the oracle if the global shrinkage parameter is of the same order as the proportion of sparisty using the asymptotic framework introduced by \cite{bogdan2011asymptotic}. Thus, it seems natural to require that any new sparse signal recovery prior should attain the oracle risk up to a multiplicative constant, and improve upon the error rates in theory as well as in practice. The generality of the Bayes risk results was conjectured by \cite{datta2013asymptotic} and proved by \cite{ghosh2013asymptotic} in a recent unpublished manuscript. \cite{ghosh2013asymptotic} proved that asymptotic Bayes optimality holds true for a general class of shrinkage priors where the local shrinkage parameter follows a distribution with a slowly-varying component bounded away from $0$ and $\infty$. This class of shrinkage priors includes many of the recently introduced priors such as the horseshoe, the normal-exponential-gamma \citep{griffin2005alternative}, the  three-parameter beta \citep{armagan2011generalized}, and the generalized double Pareto \citep{armagan2013generalized}, among others, but this class excludes the horseshoe+ prior, since its heavier tail is slowly varying but is not bounded above.  

In the light of the previous works, the purpose of our article, then, is to provide an estimator that sharpens the ability
of the Bayes estimator to extract signals from sparsity while maintaining the optimal properties of the induced decision rule.  We provide theoretical justifications by demonstrating that the proposed estimator has sharper information theoretic bounds and better MSE bounds compared to the horseshoe estimator. We illustrate
that the horseshoe+ estimator achieves greater separation of signals and noise in a standard simulation setting
and we provide a comprehensive MSE comparison with existing sparse estimators.
We develop a hierarchical model which is a natural extension of the horseshoe model
of \cite{carvalho2010horseshoe} and hence our terminology for the horseshoe+
hierarchical model.

The rest of the paper is outlined as follows. Section~\ref{sec:onetwo} motivates the class of one-group global-local mixture shrinkage priors for sparse signal estimation as a suitable alternative to the commonly used two-groups models. Section~\ref{sec:hs+} describes the
horseshoe+ estimator with a particular reference to global-local shrinkage
estimators. Section~\ref{sec:theory} provides theoretical properties of our proposed
estimator.  Our major findings can be summarized as follows: 
\ben
\item The decision rule induced by the horseshoe+ prior attains the risk of Bayes oracle under 0-1 loss up to a multiplicative constant, with the constant in Bayes risk close to the constant in oracle. We also obtain a sharper bound on the probability of type-I error compared to the horseshoe prior. 
\item The posterior mean squared error for the horseshoe+ estimator is always smaller than the posterior mean squared error of the horseshoe estimator in estimating a large signal.
\item The estimated sampling density using the horseshoe+ prior converges to the true density at a super-efficient rate when the true parameter value is zero, when the efficiency is calculated using the Kullback-Leibler (K-L) distance between the true density and the estimated sampling density. The upper bound of the risk for horseshoe+ is shown to be smaller than that of the horseshoe estimator using asymptotic properties of the prior utilizing Meijer-G functions \citep{mathai2009h}.
\een
Section~\ref{sec:sim} provides comparisons of our proposed approach with other shrinkage rules using a
standard design setting. We compare horseshoe+ with the Dirichlet-Laplace estimator \citep{bhattacharya2014dirichlet} and the horseshoe estimator \citep{carvalho2010horseshoe}, illustrating superior performance of the horseshoe+ estimator in both estimation (under squared error loss) and testing (under 0-1 loss). 
Section~\ref{sec:real} discusses the application of the proposed prior on a high-dimensional prostate cancer data set. Section \ref{sec:disc} concludes with some directions for future research.

\section{The one and two groups models}\label{sec:onetwo}

Consider the model of Section 1, i.e., $(y_i | \theta_i) \sim \Nor(\theta_i, 1)$, for $i = 1, \ldots, n$, where $\theta$ is ultra-sparse or nearly-black, in the sense that $\theta \in l_0[p_n]$. Our interest might lie in testing whether each $\theta_i$ is zero or non-zero, based on a suitably normalized test statistic or in proposing a suitable estimate $\hat \theta_i$, that has attractive properties, e.g., low mean squared error.  
The large number of parameters together with sparsity require further modeling of the data to facilitate learning via empirical Bayes or full Bayes methods. The two-groups model \citep[see, e.g.,][]{mitchell88, efron2008microarrays}, provides a natural Bayesian hierarchical framework for the sparse multiple testing problem where conditionally i.i.d. $\theta_i$ are modeled as
\beq
\theta_i | \mu = (1-\mu)\delta_{\{0\}} + \mu \Nor (0, \psi^2), \label{spikeslab}
\eeq
where $\delta_{\{0\}}$ denotes a point mass at zero and the parameter $\psi^2>0$ is the non-centrality parameter that determines the separation between the two groups. Under this setting, the marginal distribution of $y_i| \mu$ is given by
\beq
y_i | \mu \sim  (1-\mu) \Nor(0, 1) + \mu \Nor(0, 1+\psi^2). \label{twogroups}
\eeq
As can be seen from Equation \eqref{twogroups}, the two-groups model leads to a sparse estimate, i.e., it puts exact zeros in the model. The two-groups model enjoys a number of attractive theoretical properties, detailed as follows:
\begin{enumerate}
\item \cite{johnstone2004needles} showed that a thresholding-based estimator for $\theta$ under the two-groups model with an empirical Bayes estimate for $\mu$ is minimax in $\ell_2$ sense.
\item \cite{castillo2012needles} treated a full Bayes version of the problem and again found an estimate that is minimax in $\ell_2$.
\item \cite{bogdan2011asymptotic} found that the estimator under the two-groups model provides asymptotically optimal performance in testing, in the sense that its performance matches the Bayes oracle up to a constant.
\end{enumerate}
Thus, while the two-groups approach is a recognized gold-standard for Bayesian sparse signal detection and estimation, a number of arguments favor an alternative approach via the global-local shrinkage priors, also termed as the one-group model. First, in many real life applications, such as studies involving ``high-dimensional, low sample size'' gene expression data, the majority of the effect sizes are negligible, but not exactly zero, leading to an argument against exact sparsity induced by the model in Equations (\ref{spikeslab}-\ref{twogroups}). From a more pragmatic point of view, the one-group model leads to much faster computation, owing to the simple batch updating in the Gibbs sampler for the latent local shrinkage parameters. We refer the readers to \cite{carvalho2010horseshoe} for further arguments and insights. 

A useful outcome of the two-groups model is that the posterior mean $\mathbb{E}(\theta_i | y_i)$ can be written as follows: 
\beq
\mathbb{E}(\theta_i | y_i) = \omega_i \frac{\psi^2}{1+\psi^2} y_i \approx \omega_i y_i(1+o(1)) \; \text{ as } \; \psi^2 \rightarrow \infty, \label{eq:post-inc}
\eeq
where $\omega_i = P(\theta_i \neq 0 | y_i)$ is the posterior inclusion probability. 
Looking at the form of the posterior mean, one can see that it involves a ``global'' component $\psi^2/(1+\psi^2)$ that provides shrinkage towards zero for all the parameters. However, the ``local'' component $\omega_i$ allows the signal terms to escape from being too close to zero. The lack of a local shrinkage term explains why Stein-type global shrinkage estimators perform poorly in a nearly-black setting. 

The key to success in a one-group model is to design a ``global-local'' shrinkage term that gives the same form of the posterior mean as in the two-groups model. The horseshoe prior of \cite{carvalho2010horseshoe} is one such global-local shrinkage prior that has been shown to possess a number of theoretically attractive properties along with a considerably easier computational implementation compared to the two-groups model.
\begin{enumerate}
\item \cite{carvalho2010horseshoe} showed the horseshoe estimator has good information theoretic properties when the true parameter vector is sparse, in the sense that the K-L distance between the estimated and the true densities decreases at a super-efficient rate.
\item \cite{datta2013asymptotic} proved that the decision rule induced by the horseshoe estimator is asymptotically Bayes optimal for multiple testing under 0-1 loss up to a multiplicative constant.
\item \cite{van2014horseshoe} showed the horseshoe estimator is minimax in $\ell_2$ in a nearly-black case up to a constant. The constant they have been able to achieve is at least twice as large as the minimax constant of \citet{donoho1992maximum}.
\end{enumerate}
These theoretical properties, coupled with the ease of computational implementation suggests the one-group model holds considerable promise. Some other important examples of the one group model include the three-parameter beta prior \citep{armagan2011generalized}, the normal-exponential-gamma prior \citep{griffin2005alternative}, the generalized double Pareto prior \citep{armagan2013generalized}, the generalized shrinkage prior \citep{denison12} and the Dirichlet-Laplace prior \citep{bhattacharya2014dirichlet}. Below we describe the one-group horseshoe hierarchical model and then proceed to propose the horseshoe+ model that leads to considerable improvements upon the horseshoe.
\section{The horseshoe+ estimator}\label{sec:hs+}
Given normally distributed data $( y_i | \theta_i ) \sim \Nor ( \theta_i,1)$, the horseshoe hierarchical model is defined by the set of conditional distributions
\begin{align}
  ( \theta_i | \lambda_i, \tau ) &\sim \Nor \left ( 0 , \lambda_i^2 \right ),\label{eq:hs-hier}\\
  ( \lambda_i | \tau ) &\sim C^+ \left ( 0 , \tau \right ) , \nonumber
\end{align}
where $C^+$ denotes a half-Cauchy distributed scale parameter $\lambda_i$ with density 
\begin{equation}
 p(\lambda_i | \tau) = \frac{2}{\pi\tau \{1 + (\lambda_i/\tau)^2\}}, \label{eq:hs}
\end{equation}
as discussed by  \cite{gelman2006prior}. The horseshoe+ hierarchical model is defined similarly by the set of conditionals
\begin{align}
  ( \theta_i | \lambda_i, \eta_i, \tau ) &\sim \Nor \left(0 , \lambda_i^2 \right), \label{eq:hs+-hier} \\
  ( \lambda_i | \eta_i, \tau ) &\sim C^+ \left( 0 , \tau \eta_i \right),\nonumber \\
  \eta_i &\sim C^+ \left( 0 , 1 \right),\nonumber  
\end{align}
where we have introduced a further half-Cauchy mixing variable $\eta_i$. In both models, the local shrinkage random effects $\lambda_i$'s are not marginally independent after mixing over the global shrinkage parameter
$\tau$. The horseshoe+ model builds on the horseshoe by assuming that the
$\lambda_i$'s are conditionally independent given another level of local shrinkage parameters $\eta_i$'s, in addition to $\tau$.  Integrating over $\eta_i$ gives the 
the density of $\lambda_i$ as
\begin{equation}
  p( \lambda_i | \tau) = \frac{4}{\pi^2 \tau} \frac{ \log (\lambda_i/\tau)}{(\lambda_i/\tau)^2 -1}.\label{eq:lambdai}
\end{equation}
Although conceptually a natural extension, we will see that the additional $\log(\lambda_i/\tau)$ term in the numerator leads to very different properties of the proposed estimator compared to the horseshoe. There are a number of ways of dealing with the global shrinkage parameter
$\tau$. In a full Bayesian approach one can put a standard half-Cauchy prior or a $\mathrm{Uniform} (0,1)$ prior on $\tau$. Another approach is to appeal to an asymptotic argument that suggests that the
empirical Bayes estimator of $\tau$ to be set to $\hat{\tau}=p_n/n$, where $p_n$ is the
number of non-zero entries in $\theta$ \citep{van2014horseshoe}.

To further develop the distributional properties of the horseshoe+ prior we write this as a member of the class of global-local shrinkage priors with marginal prior density
$$
p(\theta_i | \tau) = \int_{0}^{\infty} p(\theta_i | \lambda_i, \tau) p(\lambda_i | \tau) d\lambda_i.
$$
Transforming to a shrinkage scale with
$\kappa_i = 1/(1 + \lambda_i^2 \tau^2)$ yields
$$ 
  p(\theta_i | \tau) = \int_0^1 p(\theta_i | \kappa_i, \tau) p(\kappa_i | \tau) d\kappa_i, 
  \; \; {\rm with} \; \; 
  p(\theta_i | \kappa_i, \tau) \sim \mathcal{N} \left(0, \frac{1-\kappa_i}{\kappa_i} \right), 
$$ 
where $\kappa_i \in [0,1]$ is a shrinkage ``weight''. 
The corresponding ``ultra-sparse'' Bayes estimator is 
\beq 
\hat{\theta}_i= \mathbb{E} ( \theta_i | y_i, \tau ) = ( 1 - \mathbb{E}(\kappa_i |y_i, \tau) ) y_i, \label{postmean}
\eeq
where we need to compute $\mathbb{E} (\kappa_i |y_i, \tau)$. 
By comparing the expression for the posterior mean for $\theta_i$ for the one-group model given by Equation \eqref{postmean} to the two-groups model given by Equation \eqref{eq:post-inc}, it is apparent that the quantity 
$\hat{\omega_i} = 1 - \mathbb{E}(\kappa_i |y_i, \tau)$ behaves as the posterior inclusion probability $P(\theta_i \neq 0 | y_i)$. This results in a natural threshold for simultaneously testing $H_{0i}: \theta_i = 0$ vs. $H_{1i}: \theta_i \neq 0$ for $i = 1,\ldots, n$. We will consider the following multiple testing procedure proposed by \cite{carvalho2010horseshoe}, and later shown to be optimal under $0$-$1$ loss by \cite{datta2013asymptotic}, for the horseshoe prior: 
\beq
\mbox{Reject } H_{0i}: \mbox{ if } 1 - \mathbb{E}(\kappa_i |y_i, \tau) > \half. \label{eq:rule}
\eeq

\subsection{Shrinkage profile}

Note that the marginal data likelihood is $p(y_i | \kappa_i, \tau) = \kappa_i^{1/2} \exp \left(-\kappa_i y_i^2/2 \right)$. Signals are identified when $\kappa_i \rightarrow 0$ and sparsity occurs when $\kappa_i \rightarrow 1$ in the posterior. We see that there are no shrinkage factors in the marginal likelihood to ``help'' identify signals in the normal model as $p(y_i | \kappa_i, \tau) \to 0$ as $\kappa_i \to 0$. This is precisely why the normal prior performs poorly for sparse settings. The horseshoe prior was designed to cancel the factor $\kappa_i^{1/2}$ and to simultaneously places prior mass at $ \kappa_i=1$ to introduce shrinkage (see \citet{carvalho2010horseshoe} for further discussion). The priors on the local shrinkage factor $\lambda_i$ and the induced prior on $\kappa_i$ for the horseshoe, the horseshoe+ and the generalized double Pareto prior are summarized in Table~\ref{tab:prior}.
\begin{table}[!ht]
\begin{center}
\caption{Priors for $\lambda_i$ and $\kappa_i$ for some global-local shrinkage rules.\\}
\label{tab:prior}
\footnotesize
\begin{tabular}{lcc}
\hline
Prior for $\theta_i$ & Prior for $\lambda_i$ & Prior for $\kappa_i$ \\ 
\hline \\
GDP & $\frac{\sqrt{2}}{(\lambda_i^2)} \int_{0}^{\infty} \exp \bigg(\sqrt{\frac{2u}{\lambda_i^2}} - u \bigg) \sqrt{u} \mathrm{d}u$  & $\frac{1}{2(1-\kappa_i)^2} \left[ \frac{\sqrt{\pi} \exp \left\{ \frac{\kappa_i}{2(1-\kappa_i)} \right\} Erfc\left\{ \sqrt{\frac{\kappa_i}{2(1-\kappa_i)}} \right\} }{\sqrt{2\kappa_i(1-\kappa_i)}} - 1 \right]$ \\[20pt]
Horseshoe & $2/ \left\{ \pi \tau (1 + (\lambda_i/\tau)^2 )\right\}$  & $\frac{\tau}{\sqrt{\kappa_i (1-\kappa_i )}} \frac{1}{(1+\kappa_i (\tau^2 -1 ) )}$ \\[10pt]
Horseshoe+ & ${4\log (\lambda_i/\tau )}/{\left\{{\pi^2 \tau}((\lambda_i/\tau)^2 -1)\right\}}$ &  $\frac{\tau}{\sqrt{\kappa_i (1-\kappa_i )}}\frac{\log \left \{ ( 1 - \kappa_i ) / \kappa_i \tau^2 \right \}}{ (1-\kappa_i (\tau^2 +1 ))}$ \\[20pt]
\hline 
\end{tabular}
\end{center}
\end{table}

\begin{center}
\begin{figure}[!t]
\centering
\includegraphics[width=10cm, height=6.4cm]{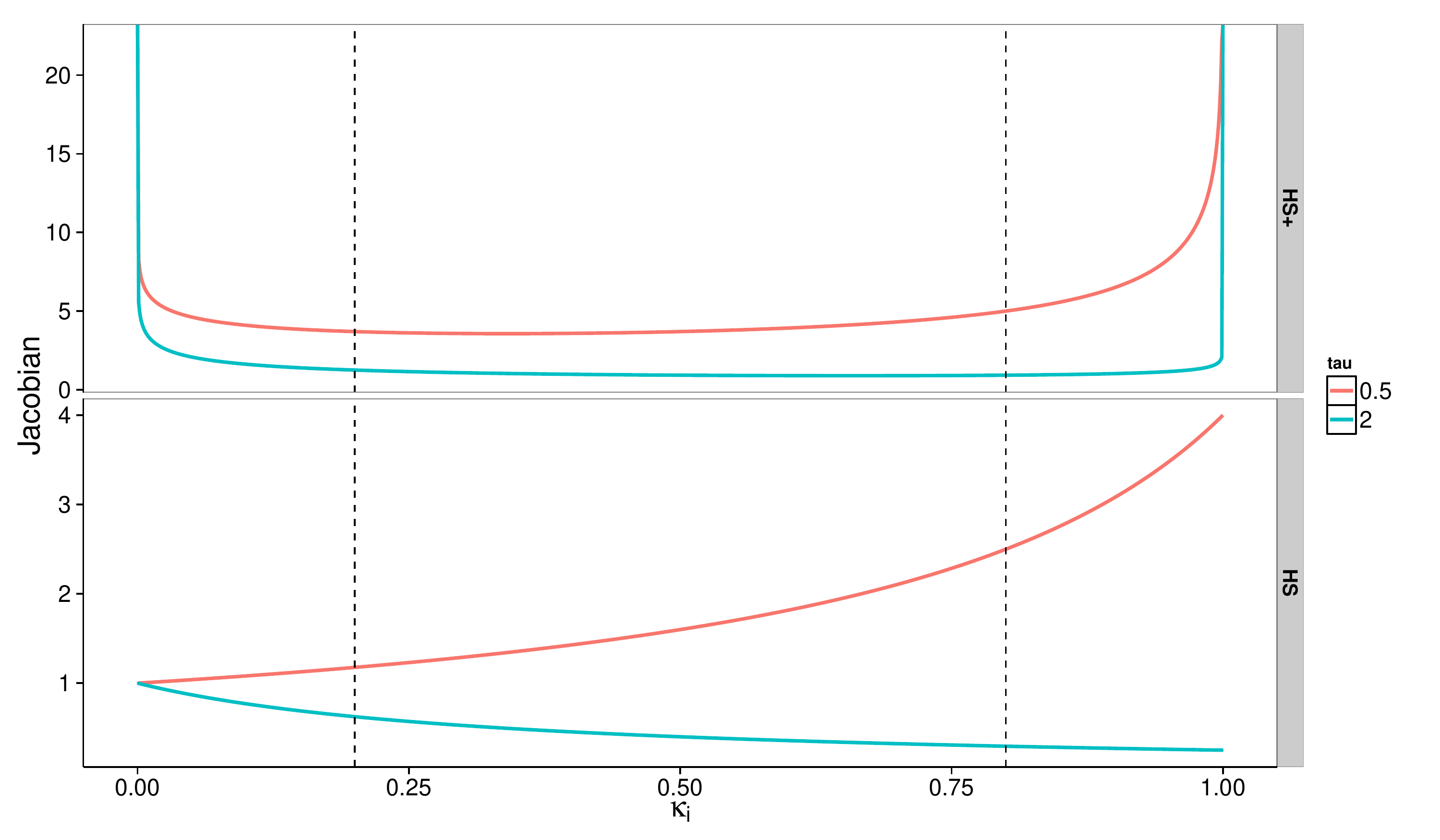}
\caption{The horseshoe+ (top) and horseshoe (bottom) prior Jacobian terms against $\kappa_i$ for $\tau=0.5$ and $2$.  The vertical lines are at $\kappa=1/(1+\tau^2)$. }
\label{fig:jacobian}
\end{figure}
\end{center}

The main difference between horseshoe+ and the others is in the extra Jacobian term introduced in the representation on the shrinkage scale. This term has a fundamentally different behavior for separating signals ($\kappa_i=0$) from the noise terms ($\kappa_i=1$). The horseshoe+ prior introduces another horseshoe $U$-shaped Jacobian factor that pushes posterior mass to the places of most interest, $\kappa_i=0,1$. This provides horseshoe+ prior with an additional power to detect signals in the ultra sparse signal case. Figure \ref{fig:jacobian} plots the Jacobians of the horseshoe and horseshoe+ priors with $\tau$ set to $0.5$ and $2$ to make the difference explicit. The horseshoe Jacobian displays unequal shrinkage behavior near the two extremities of $\kappa_i$. 

\vspace{-0.9cm}

This extra shrinkage through the Jacobian is an unique property of the horseshoe+ prior not shared by any of the other univariate shrinkage priors. To see this, note that the `sensible' priors can be expressed in terms of a slowly varying function following Theorem 1 of \cite{polson2010shrink}: 
\begin{align}
p(\lambda_i^2) & \propto (\lambda_i^2)^{(-a-1)} L(\lambda_i^2) \mbox{ for } \tau^2 = 1 \label{eq:slowly-1} \\
p(\kappa_i) & \propto (1-\kappa_i)^{(-a-1)} \kappa_i^{(a-1)} L(1/\kappa_i -1) \label{eq:slowly-2}
\end{align}
where $L(\cdot)$ is a slowly varying function with the property $L(ty)/L(y) \to 1$ as $y \to \infty$. In a recent unpublished manuscript, \cite{ghosh2013asymptotic} showed that the popular shrinkage priors like the three-parameter beta (TPB, \cite{armagan2011generalized}), which includes the popular Strawderman-Berger prior, the horseshoe prior and the normal-exponential-gamma prior, as well as the generalized double Pareto (GDP, \cite{armagan2011generalized}) prior fall into this class. Furthermore, the authors proved that the slowly varying component of \eqref{eq:slowly-1} is bounded as $\lambda_i^2 \to \infty$ for these popular shrinkage rules, i.e. $\lim_{\lambda_i^2 \to \infty} L(\lambda_i^2) \in (0,\infty)$ for priors such as TPB and GDP. This is where the horseshoe+ prior stands out from the rest, as the slowly-varying component for the prior density $p_{HS+}(\lambda_i^2)$ is unbounded as $\lambda \to \infty$, i.e. 
$$
\lim_{\lambda_i^2 \to \infty} L_{HS+}(\lambda_i^2) = \lim_{\lambda_i^2 \to \infty} \log(\lambda_i^2) \bigg( 1-\frac{1}{\lambda_i^2} \bigg)^{-1} \to \infty.
$$
Since $\kappa_i \to 0$ as $\lambda_i^2 \to \infty$, the unboundedness of  $L_{HS+}(\lambda_i^2) \equiv L(1/\kappa_i -1)$ together with \eqref{eq:slowly-2} implies that the extra shrinkage at $\lim_{\kappa_i \to 0} p(\kappa_i) \to \infty$ only holds for the Horseshoe+ prior among all shrinkage priors expressible as heavy-tailed Gaussian scale mixtures. The Jacobian term can also be interpreted on the shrinkage scale. Specifically, for $\kappa = 1/(1+\tau^2)$, we have
$$
p( \kappa_1 , \ldots , \kappa_p | \kappa , y) \propto  \prod_{i=1}^n \frac{1}{\sqrt{1-\kappa_i}} \exp\left\{ - \kappa_i \frac{y_i^2}{2} \right\}
 \frac{|\log\left((1 - \kappa_i^{-1})/(1-\kappa^{-1})\right)|}{|\kappa - \kappa_i |} .
$$
This representation shows that the horsehsoe+ prior allows differential shrinkage for $\kappa_i$ around $\kappa$ (and is continuous at $\kappa_i=\kappa$), and suggests that the global shrinkage parameter $\tau^2$ can also be interpreted as a scaling factor for the shrinkage weights $\kappa_i$.


\section{Theoretical properties of the horseshoe+ estimator}\label{sec:theory}

In this section we establish a few theoretical properties for the proposed prior and the resulting posterior, from both a decision theoretic and information theoretic viewpoint. We present our main results in the form of seven theorems. Proofs and technical details are given in Appendix \ref{app:th1}-\ref{app:mse}. 

\subsection{Marginal density for the horseshoe+ prior}\label{sec:marginal}
We start by formally establishing that the marginal prior density for horseshoe+ is unbounded at the origin.
\begin{figure}[!t]
\begin{center}
\includegraphics[width=12cm, height=7cm]{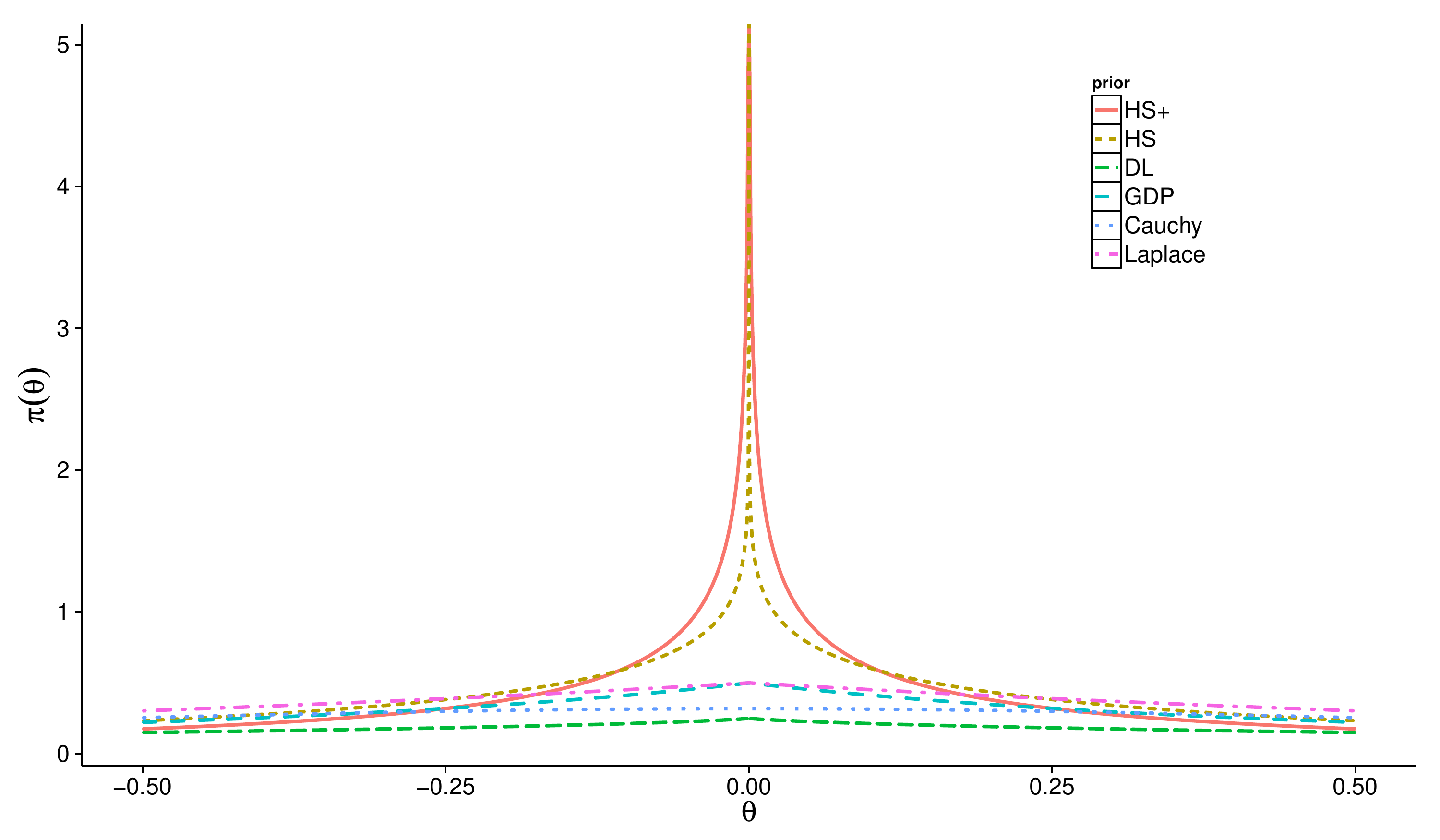}
\caption{Marginal prior densities near the origin. The legends denote the horseshoe+ (HS+),  horseshoe (HS), Dirichlet-Laplace (DL),  generalized double Pareto (GDP), Cauchy and Laplace priors.}
\label{fig:priors-origin}
\end{center}
\end{figure}
\begin{figure}[!h]
\begin{center}
\vspace{-0.5cm}
\includegraphics[width=12cm, height=7cm]{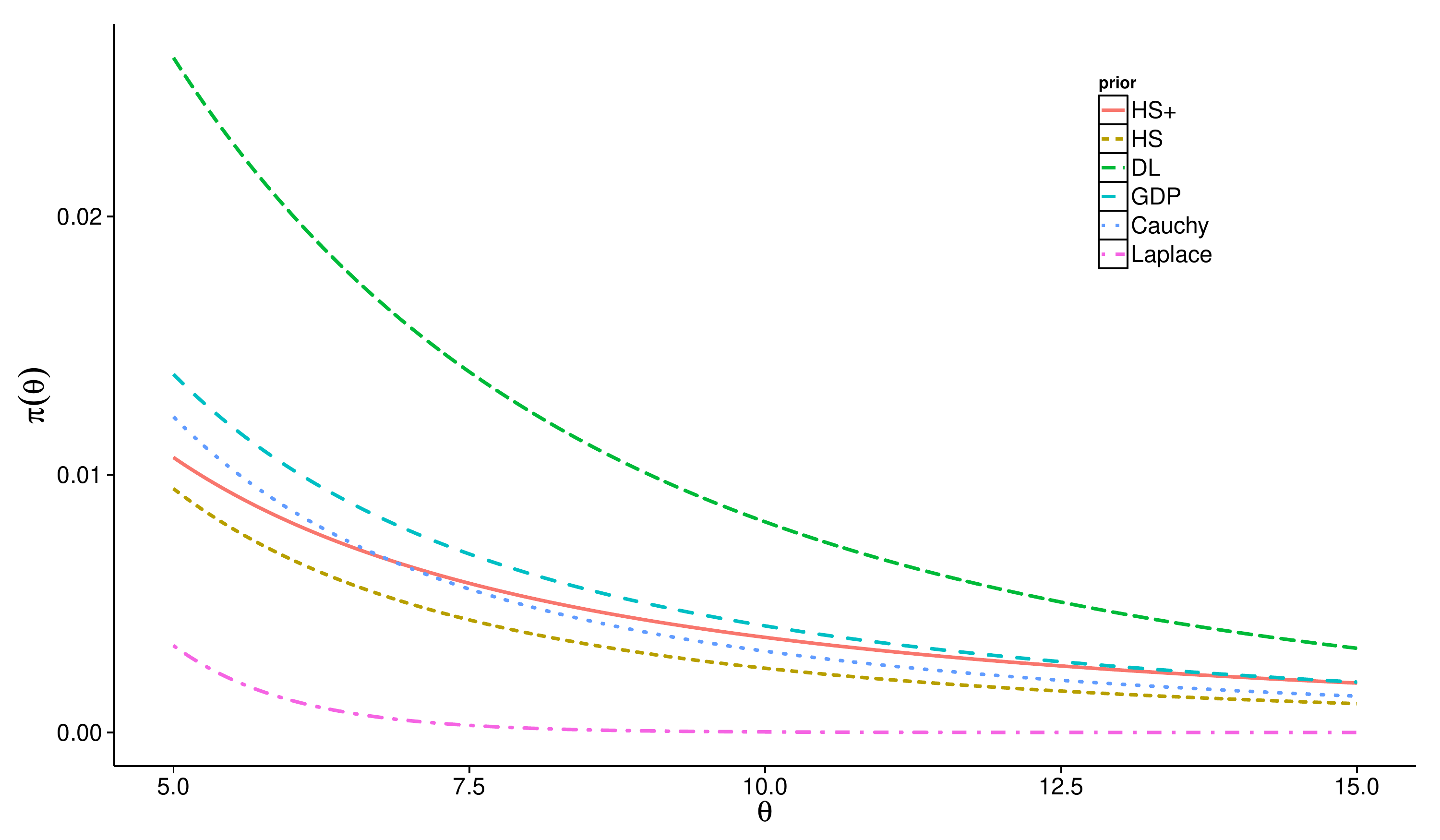}
\caption{Marginal prior densities in the tail regions. The legends denote the horseshoe+ (HS+),  horseshoe (HS), Dirichlet-Laplace (DL),  generalized double Pareto (GDP), Cauchy and Laplace priors.}
\label{fig:priors-tail}
\end{center}
\end{figure}
\begin{theorem}\label{th:1}
Assume $\tau^2=1$. Then the marginal density of the horseshoe+ prior, $p_{HS+}(\theta)$, satisfies the following properties:
\begin{enumerate}
\item 
\beq
\frac{1}{\pi^2\sqrt{2\pi}} \log \left(1+\frac{4}{\theta^2} \right) < p_{HS+}(\theta) \leq \frac{1}{\pi^2 |\theta|} \nonumber
\eeq
\item 
\beq
\lim_{|\theta| \to 0} p_{HS+}(\theta) = \infty. \nonumber
\eeq
\end{enumerate}
\end{theorem}
A proof is given in Appendix \ref{app:th1}. Figures \ref{fig:priors-origin} and \ref{fig:priors-tail} show the behavior of several global-local shrinkage priors near the origin and at the tails. The priors considered here are: horseshoe+, horseshoe \citep{carvalho2010horseshoe}, Dirichlet-Laplace \citep{bhattacharya2014dirichlet}, generalized double Pareto \citep{armagan2013generalized}, standard Cauchy, and standard Laplace (double-exponential). Note that the horseshoe+ and horseshoe densities are unbounded near the origin and more importantly, horseshoe+ puts more mass compared to the horseshoe in a small neighborhood of the origin and has heavier tails. \cite{carvalho2010horseshoe} established that a prior with unbounded density near the origin leads to super-efficiency in density estimation in a sparse signal setting. Due to Theorem~\ref{th:1}, the horseshoe+ estimator enjoys the resultant advantages as we shall show in Section \ref{sec:kl}.

\subsection{Asymptotic Bayes optimality under sparsity} \label{sec:abos}

\cite{datta2013asymptotic} proved that the Bayes risk optimality for the horseshoe prior leverages the fact that the shrinkage weight $1-\hat{\kappa}_i$ concentrates near one (uniformly in $y_i$) if the global shrinkage parameter $\tau \to 0$, and concentrates near zero if  $|y_i| \to \infty$ for any fixed $\tau$ in $(0,1)$. To attain the Bayes risk of the oracle, one additionally needs the global shrinkage parameter $\tau$ to adapt to the underlying proportion of non-zero effects $\mu_n$, i.e. $ \lim_{n \rightarrow \infty} \tau\mu_{n}^{-1} \in (0,\infty)$, where $\mu_n = \#\{\theta_i \neq 0\}/n$. It turns out that similar concentration inequalities, but with sharper bounds, hold for the posterior distribution of $\kappa_i$ under the new horseshoe+ prior. At an intuitive level, this suggests that the decision rule induced by the horseshoe+ prior will also inherit the same, if not better, optimality properties. In this section, we state two posterior concentration inequalities along with the asymptotic type-I and type-II error probabilities to establish the oracle property for horseshoe+. 
\par

Below, we briefly describe the notion of Bayes oracle in the context of multiple testing following the asymptotic framework of \cite{bogdan2011asymptotic}. Assume the two-groups model of Equations (\ref{spikeslab}-\ref{twogroups}). The optimal Bayes rule under a $0$-$1$ additive loss for testing $H_{0i}: \theta_i=0$ vs. $H_{1i}: \theta_i \neq 0$ is given by: 
\begin{align}
\mbox{Reject  } H_{0i} & \; \mbox{if} \; | y_i| > C \nonumber
\end{align}
where, 
\begin{align}
C^2 = C^2_{\psi,f} & = \frac{1+\psi^2}{\psi^2} \left ( \log (\psi^2 +1)+2 \log f \right ) \mbox{ where } f = \frac{1-\mu}{\mu}.\label{eq:c}
\end{align}
We call this rule the \textit{Bayes oracle} as the risk for this is the lower bound of $(1/n)$ times the risk for any multiple testing procedure under the two-groups model. \cite{bogdan2011asymptotic} further re-parametrized this by $u = \psi^2$ and $v = u f^2$, to obtain the following simpler form for the threshold in the oracle:
\beq
C^2 = (1+\frac{1}{u})(\log v + \log (1+\frac{1}{u})). \label{eq:oracle}
\eeq
For maintaining clarity of notations and preserving correspondence with the original work of \cite{bogdan2011asymptotic}, we use the same asymptotic framework in form of the following assumption: 
\begin{assn}\label{assnA}
The sequence of vectors $\gamma_n = (\psi_n, \mu_n)$ satisfies the following conditions:  
\begin{align*}
\mu_n & \rightarrow  0; u_n \doteq \psi^2_n \rightarrow \infty; v_n \doteq u_n f^2_n \doteq \psi^2_n \left(\frac{1-\mu_n}{\mu_n}\right)^2 \rightarrow \infty; \\
& \frac{\log v_n}{u_n} \rightarrow C \in (0, \infty) \mbox{ as }  n \rightarrow \infty.
\end{align*}
\end{assn}

\begin{remark}
The asymptotic framework provides a natural way to study the properties of the Bayes risk as the parameter vector $\gamma = (\psi, \mu)$ defining the Bayes oracle in Equation \eqref{eq:c} varies through an infinite sequence indexed by the number of tests $n$ increasing to infinity. To reduce notational complexity, we will suppress the index $n$ from $\gamma_n, \mu_n, \tau_n, \psi_n$ throughout the remainder of this section. The statements such as $\mu \rightarrow 0$ should imply that $\mu_n \rightarrow 0$ as $n \rightarrow \infty$. 
\end{remark}

\begin{remark}
It should be pointed out that the conditions are not restrictive, and are in fact minimal conditions for optimality in some sense. On one hand, the Bayes Oracle is no better than a coin toss if the limit $\psi^2/2\log(1/\mu) \to \infty$, making the test powerless, and has zero type-II error when the limit goes to zero, which could happen if one has an infinite number of replicates. The interesting cases are obtained for a finite, non-zero limit which \cite{bogdan2011asymptotic} term as ``verge of detectability'' and our results pertain to this situation. 
\end{remark}

Under Assumption \ref{assnA}, the type-I and type-II error probabilities of the Bayes oracle are given by (\cite{bogdan2011asymptotic}):
\begin{align}
t_1^{\mbox{\small BO}} & = \textrm{e}^{-C/2} \sqrt{\frac{2}{\pi v \log v}}(1+o_n), \nonumber\\ 
t_2^{\mbox{\small BO}} & = (2\Phi(\sqrt{C})-1)(1+o_n), \nonumber \\
R_{opt} & = n\left((1-\mu)t_1^{BO} + \mu t_2^{BO}\right) = n\mu(2\Phi(\sqrt{C})-1)(1+o_n), \label{optrisk}
\end{align}
where $o_n$ denotes an infinite sequence of terms, indexed by $n$ (the number of tests), converging to zero as $n \to \infty$. The last expression follows from the fact that the Bayes risk for a fixed-threshold multiple testing rule is given by $R = n\left((1-\mu)t_1 + \mu t_2\right)$ for an additive $0-1$ loss, when $t_1,t_2$ denote the type-I and type-II error probabilities respectively. A decision rule is said to attain the asymptotic Bayes optimality under sparsity (or, ABOS) if the ratio of the Bayes risk of the decision rule to the risk of the Bayes oracle (Equation~\ref{optrisk}) goes to $1$ as multiplicity $n \to \infty$. 
Now, we present the first concentration inequality on the posterior distribution of $\kappa_i$ providing the conditions under which the posterior mass of $\kappa_i$ concentrates near one. We show that an upper bound to the he posterior mass of $\kappa_i \in (0, \epsilon)$, decays as $\tau^2$. 
\begin{theorem}\label{th:3}
Suppose we have observations $y_1, \ldots, y_n$ where $y_i \sim \Nor(\theta_i, 1)$, for $i = 1,\ldots,n$, and the prior on $\theta_i$ is distributed as horseshoe+ with the hierarchical model given by (\ref{eq:hs+-hier}). Then the posterior distribution of $\kappa_i = (1+\lambda_i^2\tau^2)^{-1}$ given $y_i$ and $\tau$ satisfies the following: 
\beq
\P(\kappa_i < \epsilon | y_i, \tau) \leq e^{\frac{y_i^2}{2}} \tau^2 \epsilon(1-\epsilon)^{-2},
\eeq
for any fixed $\epsilon \in (0,1)$, and any $\tau \in (0,1)$. 
\end{theorem}
The proof is given in Appendix~\ref{app:th3}. Theorem \ref{th:3} implies that the posterior distribution of $\kappa_i$ given $\tau$ and the observation $y_i$ would converge to a point mass one if $\tau \to 0$. This leads to the following bound on the probability of type-I error rate for horseshoe+ prior, with proof given in Appendix~\ref{app:typeI}.
\begin{theorem}\label{th:typeI}
Suppose we have observations $y_1, \ldots, y_n$ from the `two-groups' model in Equation \eqref{twogroups}, and we want to test $H_{0i}: \theta_i = 0$ vs. $H_{1i}: \theta_i \neq 0$, using the decision rule of Equation \eqref{eq:rule} induced by the horseshoe+ prior. Suppose furthermore that Assumption \ref{assnA} holds for the parameter vector $(\psi,\mu)$, then the probability of type-I error for horseshoe+ decision rule is given by: 
\beq
  t_1 \leq \sqrt{\frac{2}{\pi}}\frac{\tau^2}{\sqrt{\log(1/2\tau)}}(1+o(1)).\nonumber
\eeq
\end{theorem}
\begin{remark}
It should be noted that one of the bounds (and the type-I error rate) obtained for the horseshoe+ prior are sharper than that obtained for the horseshoe prior. Theorem \ref{th:3} shows $\P_{HS+}(\kappa_i < \epsilon | y_i, \tau)= O(\tau^2)$ whereas \cite{datta2013asymptotic} obtained  $\P_{HS}(\kappa_i < \epsilon | y_i, \tau) = O(\tau)$. This relative gain will not affect the asymptotic order of the total Bayes risk derived here, but this result has interesting implications (e.g. lower false positives) nonetheless. 
\end{remark}
We now present the second concentration inequality in the other direction, with a proof in Appendix~\ref{app:th5}.
\begin{theorem}\label{th:5}
Suppose we have observations $y_1, \ldots, y_n$ where $y_i \sim \Nor(\theta_i, 1)$, for $i = 1,\ldots,n$, and the prior on $\theta_i$ is distributed as horseshoe+ with the hierarchical model given by Equation (\ref{eq:hs+-hier}). Then the posterior distribution of $\kappa_i = (1+\lambda_i^2\tau^2)^{-1}$ given $y_i$ and $\tau$ satisfies the following: 
\beq
\P(\kappa_i > \eta | y_i, \tau) \leq  e^{-\eta(1-\delta)\frac{y_i^2}{2}}\frac{1}{\tau^2}C(\eta,\delta), 
\eeq
for any fixed $\eta \in (0,1)$, any fixed $\delta \in (0,1/\eta(1+\tau^2))$ and uniformly in $y_i \in \mathbb{R}$, where $C(\eta,\delta)$ is a constant independent of $y_i$. 
\end{theorem}
A corollary of Theorem \ref{th:5} is that the posterior distribution of $\kappa_i$ given $\tau$ and  $y_i$ would converge to a point mass at zero if $|y_i| \to \infty$. 

A crucial step for proving the optimality for the horseshoe prior is the choice of the global shrinkage parameter $\tau$. \cite{datta2013asymptotic} chose $\tau$ to be of the same order as the proportion of signals $\mu$, i.e. $\tau = \tau_n = O(\mu_n)$. They also argued that the optimality of the decision rule induced by the horseshoe prior depends on how well the sparsity is captured in the hyper-parameter $\tau$. This was further supported by \cite{van2014horseshoe} who showed that the condition $\tau = O(\mu)$ is a sufficient condition for the minimaxity properties of the horseshoe estimator. Since the role of $\tau$ as a global scale parameter for the prior on local shrinkage parameters $\lambda_i$ does not change with the horseshoe+ prior, intuitively the same choice on $\tau$ would lead to the optimal type-II error rates. Under this choice of $\tau$, it follows that the type-II error for horseshoe+ decision rule has the same asymptotic order as that of the type-II error rate for the Bayes oracle. Let $C$ denote the constant in the expression for the risk of the Bayes oracle as appears in Equation \eqref{eq:oracle}. Then it follows from Theorem \ref{th:5} that the type-II error rate has the following upper bound: 

\begin{theorem}\label{th:typeII}
Suppose we have observations $y_1, \ldots, y_n$ from the `two-groups' model in Equation \eqref{twogroups}, and wish to test $H_{0i}: \theta_i = 0$ vs. $H_{1i}: \theta_i \neq 0$, using the decision rule of Equation \eqref{eq:rule}. Suppose furthermore that Assumption \ref{assnA} holds for the parameter vector $(\psi,\mu)$, and the global shrinkage parameter $\tau$ decreases to zero such that $\tau = O(\mu)$. Then for all $\eta \in (0,1)$ and $\delta \in (0,1/\eta(1+\tau^2))$, the probability of type-II error of the decision rules induced by the horseshoe+ prior is bounded above by: 
$$
t_2 \leq \left(2\Phi(\sqrt{\frac{2}{\eta(1-\delta)}} \sqrt{C}) - 1\right) (1+ o(1)).
$$
\end{theorem}
The proof is given in Appendix~\ref{app:typeII}. The proof of this theorem follows similar steps as the proof of type-II error rate for horseshoe prior in \cite{datta2013asymptotic}, where a fixed $\eta = 1/4$ and $\delta = 1/9$ were used for deriving an explicit expression.  Then it follows from Theorems \ref{th:typeI} and \ref{th:typeII} that the risk of the horseshoe+ decision rule is given by 
\begin{align*}
R_{\mbox{\small HS+}} & = n \left\{ \mu (2\Phi(\sqrt{\frac{2}{\eta(1-\delta)}} \sqrt{C}) - 1) + (1-\mu) \frac{\sqrt{2}\tau^2}{\sqrt{\pi\log(1/2\tau)}} \right\}(1+o(1)) \\
        & = n \left\{ \mu (2\Phi(\sqrt{\frac{2}{\eta(1-\delta)}} \sqrt{C}) - 1) \right\}(1+o(1)) \quad \text{as}\quad  \tau\to0.
\end{align*}
Since the risk of the Bayes oracle is $R_{\mbox{\small BO}} = n \left\{ \mu (2\Phi(\sqrt{C}) - 1) \right\}(1+o(1))$, it follows that the horseshoe+ decision rule attains the Bayes oracle up to a multiplicative
constant. 


\subsection{Kullback-Leibler risk bounds} \label{sec:kl}
\cite{carvalho2010horseshoe} proved that for horseshoe the Bayes estimate for the sampling
density, measured using the Kullback-Leibler distance between the true model
and the estimator of the density function, converges to the truth at a
super-efficient rate.  Let $\theta_0$ be the true parameter value and $f(y|\theta)$ be the sampling model. Further, let $K(q_1, q_2) = \E_{q_1} log(q_1/q_2) $ denote the K-L divergence of a density $q_2$ from $q_1$. The proof utilizes the following result by \citet{clarke90}.
\begin{proposition}\citep{clarke90}. \label{prop:clarke}
 Let $\nu_n (d\theta | y_1, \ldots, y_n)$ be the posterior distribution corresponding to some prior $\nu(d\theta)$ after observing data $y_{(n)} = (y_1, \ldots, y_n)$ according to the sampling model $f(y|\theta)$. Define the posterior predictive density $\hat q_n (y) = \int f(y|\theta) \nu_n (d\theta | y_1, \ldots, y_n)$. Assume further that $\nu(A_{\epsilon}) >0$ for all $\epsilon >0$. Then the Ces\`aro-average risk of the Bayes estimator, defined as $R_n \equiv n^{-1} \sum_{j=1}^{n} K(q_{\theta_0}, \hat q_j)$, satisfies
 \begin{equation*}
 R_n \leq \epsilon - \frac{1}{n} \log \nu(A_\epsilon),
\end{equation*}
where $\nu(A_{\epsilon})$ denotes the measure of the set $\{\theta : K(q_{\theta_0}, q_{\theta}) \leq \epsilon\}$.
\end{proposition}
Using the above proposition,  Theorem 4 of \cite{carvalho2010horseshoe} proves that for the horseshoe estimator the Ces\`aro-average risk satisfies
\begin{equation}
 R_n = O \left (\frac{1}{n} \log \left ( \frac{n}{(\log n)^b} \right) \right ), \label{eq:Rn}
\end{equation}
when the true parameter $\theta_0 = 0$. This rate is faster than any
prior without a pole at zero. It is super-efficient, in the sense that the risk is lower than that of the MLE, which has the rate $O(\log n/n)$. The same result holds for the horseshoe+ estimator due to its infinite mass near zero (by Theorem \ref{th:1}). However, we demonstrate that the horseshoe+ prior in fact has a better rate of convergence than the horseshoe prior. Our result is based on the following theorem.

\begin{theorem}\label{th:int}
Let $p^0_{HS+} (\theta)$ and $p^0_{HS} (\theta)$ denote the marginal densities of the horseshoe+ and horseshoe priors at the origin. Then we have
$$
\int_{0}^{\frac{1}{\sqrt{n}}} p^0_{HS+} (\theta) d \theta = \frac{1}{\sqrt{2} \pi^{5/2} \sqrt{n}} \left( \frac{\log^2(n)}{4} + \left( 1- \frac{\gamma}{2} + \frac{\log(4)}{4}\right) \log(n) + O(1) \right),
$$
where $\gamma$ is the Euler-Mascheroni constant and
$$
\int_{0}^{\frac{1}{\sqrt{n}}} p^0_{HS} (\theta) d \theta = \frac{1}{\sqrt{2} \pi^{3/2} \sqrt{n}} \left( \frac{\log(n)}{2}  + O(1) \right). 
$$
\end{theorem}

The proof is given in Appendix~\ref{app:meijer}. Due to the extra $\log(n)$ factor, the horseshoe+ prior places more mass around a neighborhood of the origin compared to the horseshoe prior. 
Thus, when $\theta_0=0$, setting $\epsilon = 1/n$ gives after some algebra from Proposition~\ref{prop:clarke} that
$$
R_n (HS)  \leq \frac{\log n}{2n} + \frac{1}{n} - \frac{\log \log n}{n} + \mathrm{const.}
$$
and
$$
R_n (HS+)  \leq \frac{\log n}{2n} + \frac{1}{n} - \frac{2\log \log n}{n}  + \mathrm{const.}
$$
Therefore, the multiplier of the $\log \log n$ term improves for horseshoe+.
\subsection{Mean squared error} \label{sec:mse}
It is well known that if $p(|y_i -\theta_i |)$ is the standard normal density and $p(\theta_i)$ is a zero mean scale mixture of normals, with the scale parameter $\lambda^2$ following a proper prior law, the posterior moments of $\theta_i$ admits the following representations, also known as ``Tweedie's formula'' \citep{efron2011tweedie}: 
\begin{align}
\mathbb{E}(\theta_i | y_i) & = y_i + \frac{d}{dy_i} \log m(y_i), \label{tweedie1}\\
\mathbb{V} (\theta | y_i) & = 1 + \frac{d^2}{dy_i^2} \log m(y_i)\label{tweedie2},
\end{align}
where $m(y_i)$ is the marginal for $y_i$ (see for example \citet{pericchi1992exact} and \citet{carvalho2010horseshoe}).
Furthermore, we can use properties of slowly varying functions to show that if the prior on $\theta_i$ can be written as a normal scale mixture with a ``slowly-varying" prior on the scale parameter, the marginal inherits the slowly varying property. For priors with a polynomially heavy tail it can also be shown that the resulting posterior mean is asymptotically robust, in that the difference $|\mathbb{E}(\theta_i |y_i, \tau) - y_i|$ vanishes for large $|y_i|$ while $\tau$ is fixed. 

Heavy-tailed distributions are often characterized by the notion of regular variation. The following definition is due to Karamata (see \cite{mikosch1999regular} or \cite{bingham1989regular} for a detailed discussion). 
\begin{Def}
A positive, measurable function $L$ is called regularly varying at infinity with index $\alpha$ if it is defined on the interval $[x_0, \infty)$ for some $x_0$ and 
$$
\lim_{x \rightarrow +\infty} \frac{L(tx)}{L(x)} = t^{\alpha} \quad \text{for all} \quad t>0 
$$
$L(\cdot)$ is called a slowly varying function at infinity if $\alpha = 0$. 
\end{Def}
Using the above definition, we state the following result from Theorem 6.1 of \cite{barndorff1982normal}.
\begin{proposition}\citep{barndorff1982normal}.\label{prop:barndorff}
Consider the Gaussian scale mixture $y|\lambda^2 \sim \mathrm{Normal} (0, \lambda^2)$ and suppose the prior density of $\lambda^2$ is given by $f(\lambda^2) = (\lambda^2)^{\alpha-1} L(\lambda^2)$ as $\lambda^2 \to \infty$, where $L(\cdot)$ is a slowly varying function. Then the marginal $m(y)$ after integrating out $\lambda^2$ has the property that $m(y) \propto |y|^{2\alpha-1} L(y^2)$ as $|y| \to \infty$.
\end{proposition}
Let $m_{HS+}(y_i)$ and $m_{HS}(y_i)$ denote the marginals under the horseshoe+ and horseshoe priors respectively. Proposition~\ref{prop:barndorff} immediately shows that we have that $m_{HS+}(y_i) = m_{HS}(y_i) \log(|y_i|)(1+o(1))$ as $|y_i|\to \infty$, since the only difference between the horseshoe and horseshoe+ mixing densities is the additional slowly varying $(\log \lambda_i)$ term in the scale mixing density for the horseshoe+ prior. 
In particular, as $|y_i| \to \infty$, we have
\begin{align*}
m_{HS+}(y_i) & = m_{HS}(y_i) \times \log(|y_i|) \times \frac{y_i^2-1}{y_i^2+1} \times {\mathrm{constant}}.
\end{align*}
where ``constant'' denotes the collection of all terms that does not involve $y_i$. Thus,
\begin{align}
\frac{d}{dy_i} \log m_{HS+}(y_i) & = \frac{d}{dy_i} \log m_{HS}(y_i) + \frac{1}{|y_i| \log |y_i|} - \underbrace{\frac{4y_i^2}{y_i^4-1}}_{O(1/y_i^2)}, \label{hseq1}
\end{align}
and, 
\begin{align}
\frac{d^2}{dy_i^2} \log m_{HS+}(y_i) & = \frac{d^2}{dy_i^2} \log m_{HS}(y_i) - \frac{1+ \log y_i}{(y_i\log y_i)^2} + O\left(1/y_i^3\right). \label{hseq2}
\end{align}
Using Equations (\ref{tweedie1}) and (\ref{tweedie2}), in combination with Equations \eqref{hseq1} and \eqref{hseq2}, allows one to relate the bias and variance, and hence the MSE, for the horseshoe and the horseshoe+ estimators. We have the following result:
%
\begin{theorem}\label{th:mse}
Suppose $p(|y_i -\theta_i|)$ is the standard normal density, and $p_{HS}(\theta_i)$ and $p_{HS+}(\theta_i)$ denote the horseshoe and horseshoe+ prior densities on $\theta_i$, leading to the posterior mean squared errors $\MSE_{HS} (\theta_i|y_i) $ and $\MSE_{HS+} (\theta_i|y_i)$ respectively. Then, for large values of $|y_i|$, we have,
\begin{align*}
\MSE_{HS+} (\theta_i|y_i) &= \MSE_{HS} (\theta_i|y_i) - \frac{1}{y_i^2\log |y_i|} + O\left(\frac{1}{y_i^3}\right).
\end{align*}
\end{theorem}
The proof is given in Appendix~\ref{app:mse}. This theorem establishes that the horseshoe+ estimator has asymptotically lower MSE compared to the horseshoe estimator when $|y_i|$ is large, due to the extra $(\log |y_i|)$ factor in the marginal, which in turn is due to the extra $(\log \lambda_i)$ term in the prior mixing density.


\section{Numerical examples}\label{sec:sim}

\subsection{Sum of squared error about the posterior median}

We follow the simulation setting described in  \cite{bhattacharya2014dirichlet}.  We simulate data $y_i | \theta_i \sim \mathrm{Normal} (\theta_i, 1)$ for
~$i = 1,\dots,n$, where $\theta_i=A$ in fraction $q$ of its components with the magnitude of $A=7, 8$ and $\theta_i=0$ in the remaining components. We report simulation results for $n=200$ in Table \ref{tab:rep-mse-200}. Each configuration is replicated $100$ times and the average sum of squared error
about the posterior median is reported. 

\begin{table}[!h]
  \caption{Average SSE about the posterior median for $n=200$ for the competing priors. The averages are computed over $100$ replicates. The lowest SSE for each setting (in rows) is in bold.}
  \label{tab:rep-mse-200}
   \footnotesize
   \begin{center}
      \begin{tabular}{llccccc}
\toprule
q & A &  D-L & HS Cauchy & HS+ Cauchy & HS Unif & \multicolumn{1}{c}{HS+ Unif} \\ 
\toprule
0.05 & 7  & $\phantom{0}26.86$ & $\phantom{0}\textbf{15.95}$ & $18.58$ & $\phantom{0}17.11$ & $\phantom{0}18.08$ \\
 & 8  & $\phantom{0}22.49$ & $\phantom{0}\textbf{14.47}$ & $15.97$ & $\phantom{0}15.26$ & $\phantom{0}17.42$ \\
 \hline
0.1 & 7  & $\phantom{0}43.76$ & $\phantom{0}33.92$ & $\textbf{31.65}$ & $\phantom{0}35.13$ & $\phantom{0}33.51$ \\
 & 8  & $\phantom{0}43.81$ & $\phantom{0}32.28$ & $\textbf{29.77}$ & $\phantom{0}33.67$ & $\phantom{0}32.23$ \\
 \hline
0.2 & 7  & $\phantom{0}78.11$ & $\phantom{0}69.29$ & $\textbf{59.26}$ & $\phantom{0}83.61$ & $\phantom{0}59.92$ \\
 & 8  & $\phantom{0}82.64$ & $\phantom{0}70.72$ & $\textbf{62.64}$ & $118.52$ & $\phantom{0}63.69$ \\
 \hline
0.3 & 7  & $103.46$ & $104.33$ & $\textbf{86.77}$ & $322.93$ & $100.26$ \\
 & 8  & $121.04$ & $108.12$ & $\textbf{93.21}$ & $373.71$ & $220.16$ \\
\hline 
\end{tabular}
      \end{center}
\end{table}

We compare the proposed horseshoe+ prior with two competitors: the horseshoe prior of \citet{carvalho2010horseshoe} and the Dirichlet-Laplace (D-L) prior of \cite{bhattacharya2014dirichlet}. To deal with the global shrinkage parameter $\tau$ for the horseshoe and the horseshoe+ priors, we try two scenarios: (a)  $\tau \sim \operatorname{C}^+(0,{1}/{n})$ and  (b) $\tau \sim \operatorname{Unif}[0,1]$. 
For posterior sampling, we use the Stan software package \citep{stan-software:2014} to draw $10,000$
samples in each case, half of which are treated as burn-in and discarded. We monitored MCMC convergence and found no evidence of mixing problems. The D-L prior is implemented in its hierarchal normal-exponential form,
and the horseshoe and horseshoe+ priors by the hierarchical model in Equations (\ref{eq:hs-hier}) and 
(\ref{eq:hs+-hier}) respectively.
\begin{figure}[!t]
\begin{center}
\includegraphics[width=12cm, height=7.1cm]{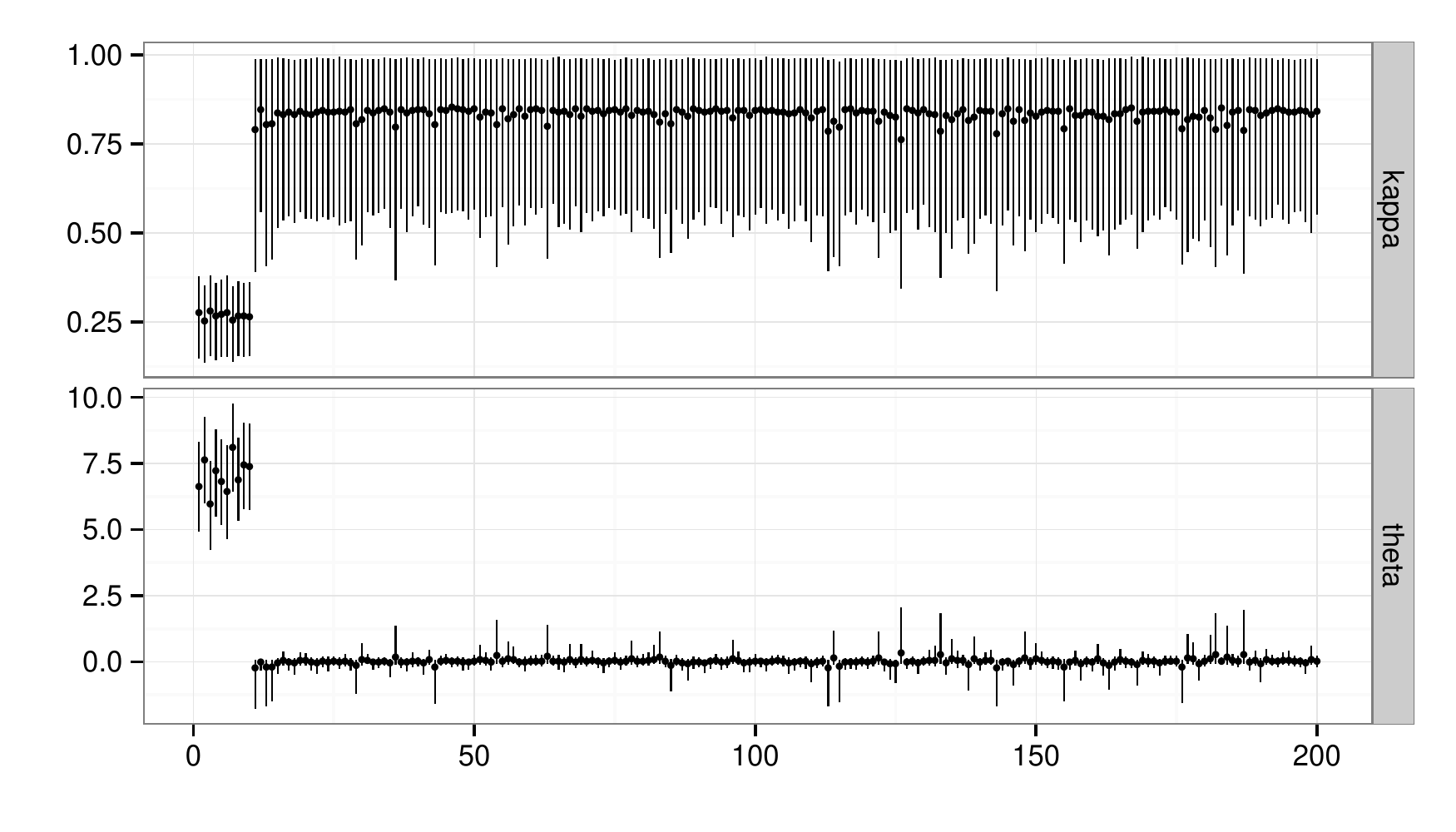}
\caption{Estimated $\kappa_i$ and $\theta_i$ for horseshoe for $n=200$ with first 10 true $\theta_i$ equal to 7 and rest true values set to 0. Dots are posterior means and solid lines are the middle 95\% posterior credible intervals. We used  $\tau \sim \operatorname{Unif}[0,1]$.}
\label{fig:hs-kappas}
\end{center}
\end{figure}
\begin{figure}[!h]
\begin{center}
\vspace{-0.2in}
\includegraphics[width=12cm, height=7.1cm]{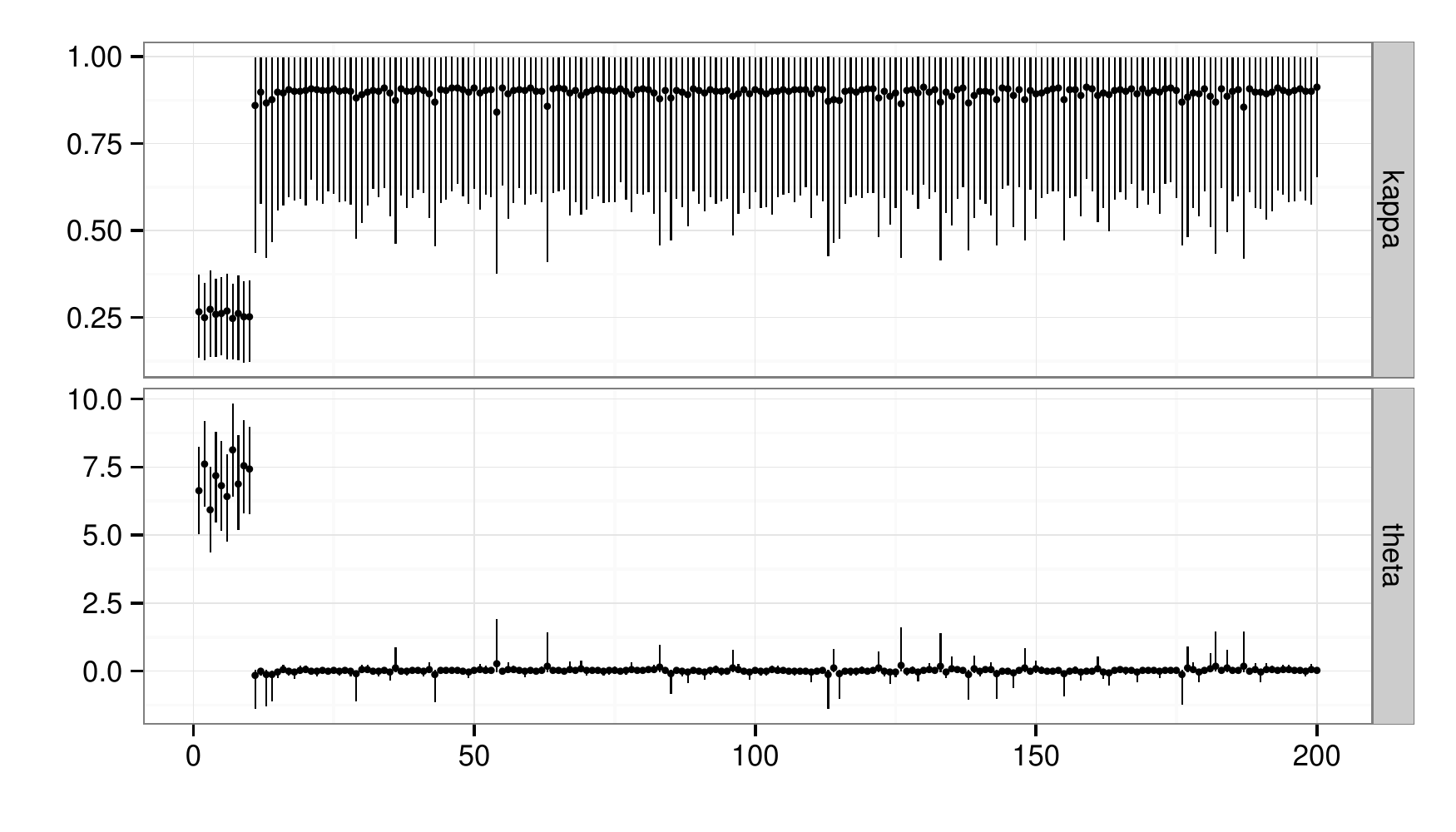}
\caption{Estimated $\kappa_i$ and $\theta_i$ for horseshoe+ for $n=200$ with first 10 true $\theta_i$ equal to 7 and rest true values set to 0. Dots are posterior means and solid lines are the middle 95\% posterior credible intervals. We used  $\tau \sim \operatorname{Unif}[0,1]$.}
\label{fig:hs+-kappas}
\end{center}
\end{figure}
\par In Table~\ref{tab:rep-mse-200}, the estimator with the lowest average SSE is in bold in each simulation setting (in rows). The horseshoe+ prior with the half-Cauchy prior on $\tau$ has the lowest SSE in all but two cases, in which the horseshoe prior performs the best. The $\operatorname{C}^+(0,{1}/{n})$ prior on $\tau$ results in better performance over a $\mathrm{Uniform}(0,1)$ prior for both horseshoe and horseshoe+ since the former puts more mass in a neighborhood close to zero, helping $\tau$ adapt to the sparsity level of the data.

\par To make the difference between the horseshoe and the horseshoe+ estimates clear, we plot  $\E(\kappa_i | y_i)$ and $\E(\theta_i |y_i)$ for $i=1, \ldots, n,$ for  horseshoe  in Figure~\ref{fig:hs-kappas} and for  horseshoe+ in Figure~\ref{fig:hs+-kappas}. In both cases, the prior on $\tau$ is $\operatorname{Unif}[0,1]$. We used $n=200$ and simulated $y_i$ with 10 components with a mean equal to 7 and the rest with mean 0.  Without loss of generality, the components (true values and estimates) with true non-zero means are plotted as the first 10 data points and those with true zero means are plotted afterwards. The posterior means are shown as dots and the middle 95\% posterior credible intervals by solid lines. By comparing the estimates, it is clear that horseshoe+ does a much better job compared to horseshoe in terms of shrinking the noise terms to zero (estimated $\hat \kappa_i$ closer to 1 or equivalently, estimated $\hat \theta_i$ closer to zero).
%
%
%
%
%
%
\subsection{Misclassification probabilities}\label{misclassprob}
We compared the performance of the multiple testing rule induced by the horseshoe+ prior with two other global-local shrinkage priors:  the horseshoe prior of \citet{carvalho2010horseshoe} and the Dirichlet-Laplace prior of \cite{bhattacharya2014dirichlet} in terms of the misclassification probability (MP). We use the misclassification probability as a criteria for our experiment as it is equal to the Bayes risk under a $0$-$1$ additive loss for data generated by a two-groups model. We follow the same experimental set up in \cite{bogdan2008comparison}, replicated in \cite{datta2013asymptotic}, where the Bayes oracle (BO) acts as the lower bound and the $MP = \mu$ line as the upper bound, where $\mu$ is the proportion of signals. We simulated data of size $n=200,\psi_n = \sqrt{2\log n}=3.26$. Our data generation scheme follows the conditions provided by \cite{bogdan2011asymptotic}, which guarantees the optimality of the Benjamini-Hochberg procedure to use it as another practical lower bound along with the Bayes Oracle. 

\par
\begin{figure}[!t]
\begin{center}
\includegraphics[width=14cm, height=7.5cm]{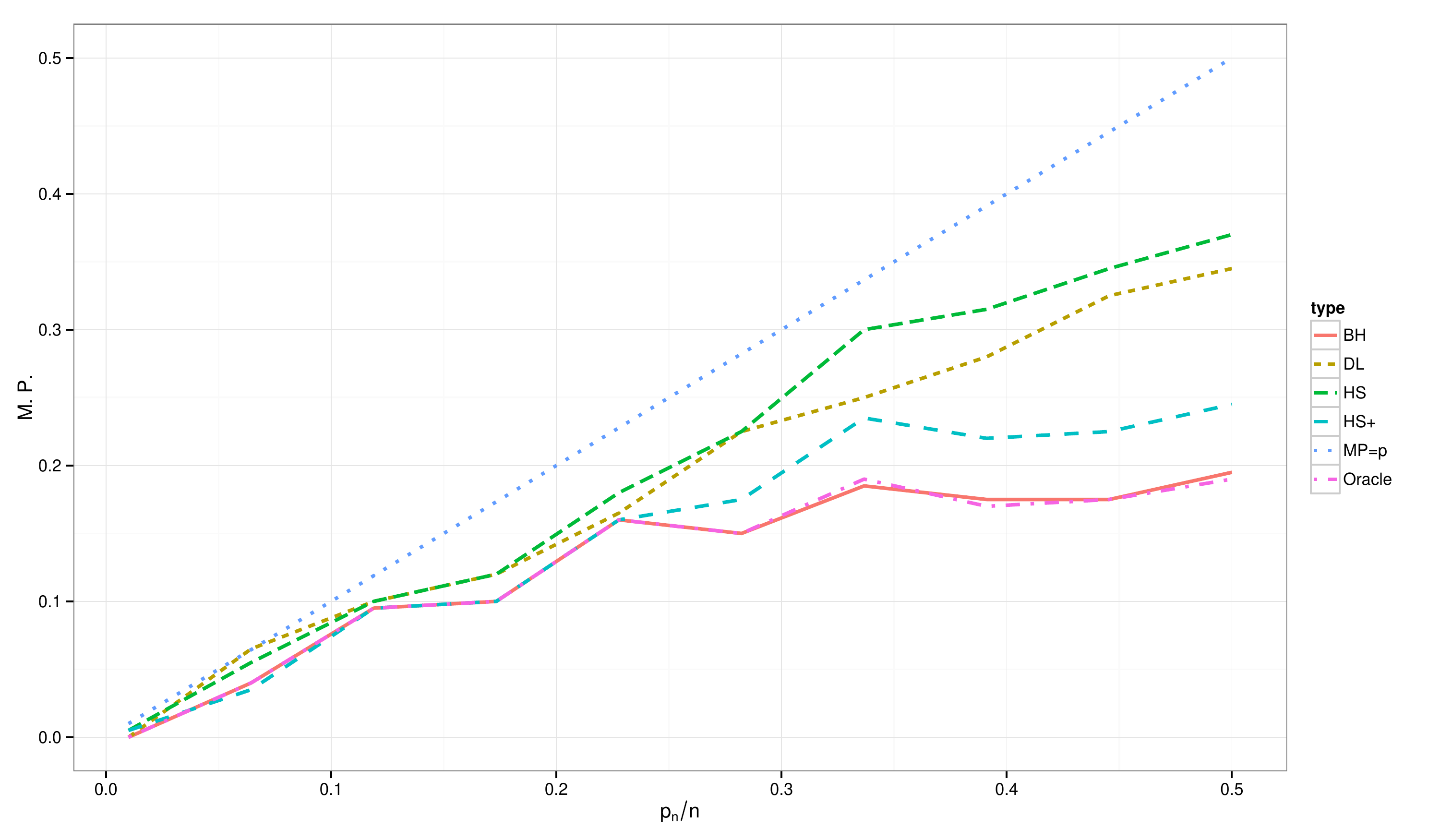}
\caption{Misclassification probability plots for the horseshoe+, horsesshoe, and the Dirichlet-Laplace ($DL_{1/n}$) shrinkage priors, Benjamini-Hochberg and the Bayes oracle for $\mu \in (0.1, 0.5)$. }
\label{fig:erplot}
\end{center}
\end{figure}

 Figure \ref{fig:erplot} shows the misclassification probabilities (henceforth abbreviated as MP) for different shrinkage priors considered for ten equispaced values of $\mu \in [0.01,0.5]$ along with the oracle and the straight line (MP = $\mu$). Figure \ref{fig:erplot} shows that the misclassification probability for the horseshoe+ prior is very close to that of the Bayes oracle for a wide range of values of $\mu$, and departs a little for values higher than $0.2$. Furthermore, the horseshoe+ decision rule leads to a superior performance compared both the horseshoe and the Dirichlet-Laplace prior. We have also plotted the MP for the Benjamini-Hochberg rule, for $\alpha = 1/\log{n} = 0.1887$, along with the one-group shrinkage priors. Under this setting, the Benjamini-Hochberg rule achieves the same MP as the oracle. This is in concordance with the theoretical results for optimality of BH in \cite{bogdan2011asymptotic}. 
\par We used the full Bayes estimates for the hyperparameters for both the horseshoe prior and the double exponential prior. For estimating $\tau$, we assumed standard half-Cauchy prior on $\tau$ for deriving the full conditionals using a Gibbs sampler. As pointed out by \cite{carvalho2009handling} and \citet{scott2006exploration}, the fully Bayesian approach for estimating $\tau$ has a few adavantages over its alternatives, viz. empirical Bayes and cross-validation. In the extremely sparse case, the empirical Bayes estimate of $\tau$ might collapse to $0$ \citep{scott2010bayes, bogdan2008comparison}. Cross-validation, though free of this problem, uses plug-in estimates for the signal-to-noise ratio. 
\cite{carvalho2009handling} argue that the plug-in estimates are not necessarily wrong, but caution should be exercised while using them for extremely sparse problems.

\section{Application on a prostate cancer data set} \label{sec:real}
We illustrate the performance of the horseshoe+ prior for the benchmark \textit{prostate cancer data}, introduced by \cite{singh2002gene} and made popular by \citet{efron2008microarrays,efron2010large,efron2010future},  among others. The \textit{prostate cancer data} has gene expression values for $n = 6,033$ genes for $m=102$ subjects, with $m_1 = 50$ normal controls and $m_2 = 52$ prostate cancer patients. The goal is to identify ``genes" that are differentially expressed between controls and the cancer patients. To analyze this data further, the test statistic values are calculated for each of the $6,033$ genes by first calculating a two-sample $t$-statistic, say $t_i, i = 1, 2, \ldots, n=6,033$ for each of the genes and then applying the inverse Normal CDF transformation to obtain $y_i = \Phi^{-1}\left( F_{t_{100 d.f.}}(t_i) \right)$. The $y_i$-values can be modeled as independent Gaussian variables with mean $\theta_i$'s, i.e. $y_i \sim \theta_i + \epsilon_i$ to cast this problem as a high-dimensional normal means inference problem. The corresponding multiple testing problem would be to simultaneously test the hypotheses $H_{0i}: \theta_i = 0,$ for $ i = 1, \ldots, n$. Under the global null hypothesis of no `differentially expressed' genes, one should expect the histogram of the test statistics to follow a $\Nor(0,1)$ density curve but the histogram shows a heavier tail, suggesting the presence of a few regulatory genes. \par 

For a proper appraisal of the extra shrinkage by the horseshoe+ prior at the tails compared to the horseshoe prior, we do the following experiment: We consider the top $10$ genes selected by \cite{efron2010future} and their effect sizes estimated by a two-groups normal hierarchical model. We apply both the horseshoe and the horseshoe+ prior to the 6,033 test statistics, and compare the `effect-size' estimates $\hat{\theta_i}$ for these genes. One would expect that the horseshoe+ prior would shrink these ``top" genes even less than the horseshoe prior and as a result the posterior mean $\hat{\theta}_i = (1-\mathbb{E}(\kappa_i | y_i, \tau))y_i$ would be closer to the observed test statistics $y_i$. \par

\begin{figure}[!t]%
\centering
\includegraphics[height=3in,width=0.75\columnwidth]{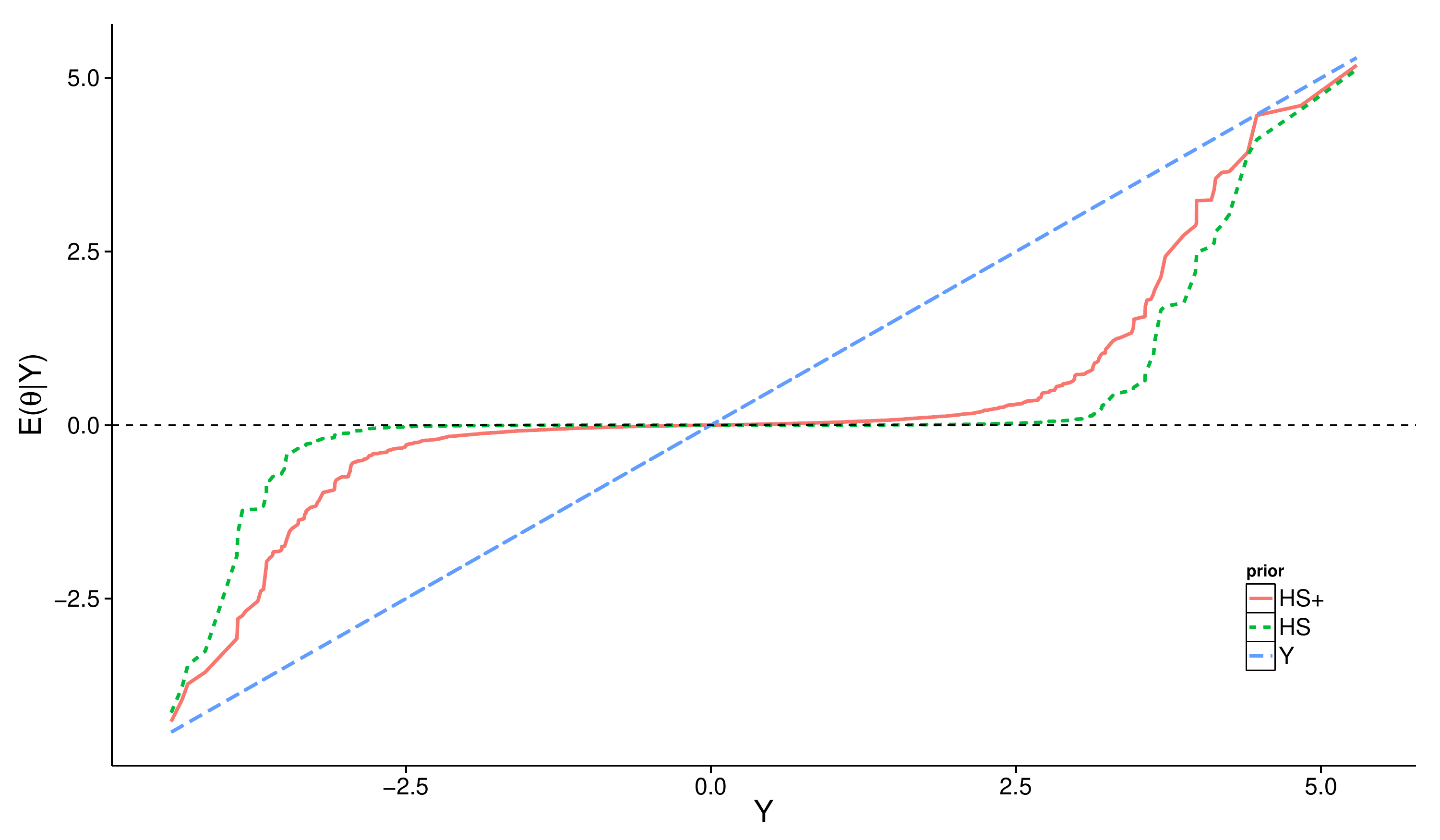}%
\caption{Posterior mean $\mathbb{E}(\theta_i | y_i)$ against $y_i$ for $6,033$ genes for the horseshoe and horseshoe+ priors applied to the prostate cancer example.}%
\label{fig:top10}%
\end{figure}

Table \ref{tab:top10} shows the top 10 genes selected by \cite{efron2010future}, and the effect size estimates by the horseshoe and the horseshoe+ priors. For both the horseshoe and horseshoe+ prior, we implemented a Gibbs sampler with 15,000 draws with a burn-in period of 3,000 draws. The benefits of a heavier tail become apparent from this table as in 9 out of the top 10 genes, the horseshoe+ estimates are closer to the observed test statistics compared to the horseshoe estimates. One might naturally wonder about the performance of the two competing Bayesian models for the ``uninteresting" genes, and it turns out that both the priors have equal strength in squelching the noisy test statistics to zero. Figure \ref{fig:top10} shows the posterior mean for the two priors 
against the observed test statistics. It can be clearly seen that all the procedures show good shrinkage properties near zero, and the only difference comes from the performance near tails, or \textit{robustness to large signals}. This is also reflected in the value of the estimated mean squared prediction error calculated as $MSE = (1/n) \sum_{i=1}^{n} (\hat{\theta}_i - y_i)^2$. The values of the mean squared prediction error for the horseshoe+ and the horseshoe prior are $1.189$ and $1.045$ illustrating the superiority of horseshoe+ prior over the horseshoe prior.

\begin{table}[t!]
  \centering
	\footnotesize
		 \caption{The test statistics ($y$-values) and the effect-size estimates for the top 10 genes selected by \cite{efron2010future} by the horseshoe, horseshoe+ models, and Efron's two-groups model estimates.}
		   \label{tab:top10}%
    \begin{tabular}{rrrrr}
    \toprule
    Gene & $y$-value & $\hat{\theta}^{HS+}_i$ & $\hat{\theta}^{HS}_i$ & $\hat{\theta}^{\mathrm{Efron}}_i$ \\
    \toprule
    610   & 5.29  & 5.20  & 5.12  & 4.11 \\
    1720  & 4.83  & 4.77  & 4.54  & 3.65 \\
    332   & 4.47  & 3.24  & 4.11  & 3.24 \\
    364   & -4.42 & -4.43 & -4.14 & -3.57 \\
    914   & 4.40  & 4.40  & 3.89  & 3.16 \\
    3940  & -4.33 & -3.78 & -3.77 & -3.52 \\
    4546  & -4.29 & -3.88 & -3.46 & -3.47 \\
    1068  & 4.25  & 3.71  & 3.03  & 2.99 \\
    579   & 4.19  & 3.99  & 2.88  & 2.92 \\
    4331  & -4.14 & -3.48 & -3.26 & -3.30 \\
    \hline
    \end{tabular}
\end{table}%

\section{Discussion}\label{sec:disc}

We have provided a default Bayesian shrinkage estimator for extracting signals from a sparse parameter vector. The proposed prior is called horseshoe+ prior as it renders itself to a natural extension of the horseshoe prior and provides substantial improvement for the ``nearly-black" or ``ultra-sparse" situation. In particular, the heavier tails of the horseshoe+ prior leads to an increasing ability of separating the signals, and the larger mass near origin leads to better handling of sparsity and a higher order of super-efficiency for the risk in density estimation. We have examined this new prior both theoretically and empirically by considering the estimation accuracy for a high-dimensional parameter vector as well as the error rates for the multiple testing rule induced by applying a threshold rule to the pseudo posterior inclusion probabilities. \par

Our asymptotic results demonstrate that the horseshoe+ estimator achieves a lower MSE and the horseshoe+ decision rule attains the Bayes oracle in testing up to $O(1)$ with a sharper bound on the type-I error rate compared to horseshoe. While we have not discussed the asymptotic minimaxity properties of the horseshoe+ estimator in the $\ell_2$ sense, we conjecture that the asymptotic minimaxity will continue to hold, likely with a sharper bound on the constant term compared to \cite{van2014horseshoe}. The sharpening effect of the horseshoe+ prior can be attributed to the extra shrinkage gained by having a U-shaped Jacobian over a lopsided one, in addition to the U-shaped prior induced on $\kappa_i$. It is also worth noting that asymptotic minimaxity results for a class of priors with  $p(\lambda_i^2) \propto (\lambda_i^2)^{-a-1} L(\lambda_i^2)$ where $L(\cdot)$ is a bounded, slowly varying function has been established by \citet{ghosh14}. For horseshoe+, $L(\lambda_i) =\lambda_i^2 \log \lambda_i /(\lambda_i^2 -1)$ is slowly varying but not bounded. It is an interesting exercise to see if the framework of \citet{ghosh14} can be adapted for horseshoe+.


In the recent past, there have been a few shrinkage priors that we collectively call the `global-local' shrinkage priors following \cite{polson2010shrink}. These priors include the generalized double Pareto \citep{armagan2013generalized}, the normal-exponential-gamma \citep{griffin2005alternative}, the three parameter beta \citep{armagan2011generalized}, and the Dirichlet-Laplace \citep{bhattacharya2014dirichlet}, among others. These priors exhibit similar shrinkage properties as the horseshoe prior in that they simultaneously squelch the noise to zero and recover the signals. Though these priors lead to competitive performances in the sparse signal recovery problem, they also have unique, distinguishable characteristics. For example, the generalized double Pareto prior leads to a closed form prior density of $\theta$ and induces a sparsity favoring penalty in regularized least squares, while the Dirichlet-Laplace prior models the joint distribution of $\theta$ under the two-groups model via the joint distribution of the shrinkage parameters. The behavior of the marginal prior densities on $\theta$ can be seen from Figures \ref{fig:priors-origin} and \ref{fig:priors-tail}, and our simulation results suggest improvements for both estimation and testing, but it is an open question whether a set of necessary conditions can be imposed on the class of global-local shrinkage priors that guarantees certain desirable properties. 


A key insight we gain from the success of the family of the global-local shrinkage priors is that the global shrinkage parameter plays a vital role in controlling the behavior of the posterior. Specifically, the global shrinkage parameter in the horseshoe prior needs to be of the order of the proportion of non-null effects to ensure asymptotic minimaxity in estimation \citep{van2014horseshoe} as well as the optimality of the induced decision rule in testing \citep{datta2013asymptotic}. We have proved that the same condition also guarantees the optimal performance for the horseshoe+prior in testing.

Finally, the horseshoe+ prior can be further extended by modeling the local shrinkage parameter $\lambda_i$ as a higher order product of independent half-Cauchy random variables, leading to an even heavier tail and larger spike at zero. The moments and densities of the Cauchy product $C_1C_2 \ldots C_k$ are given in \cite{bourgade2007euler}. The density $\Psi_k(\cdot)$ of the $k$-product $C_1C_2 \ldots C_k$ for the even and the odd cases are as follows:
\begin{align*}
\Psi_{2i+1}(x) & = \frac{2^{2i}}{\pi(2i)!}\left( \prod_{j=1}^{i} ((j-\half)^2+ \frac{(\log|x|)^2}{\pi^2}) \right)\frac{1}{1+x^2}, \\
\Psi_{2i}(x) & = \frac{2^{2i-1}}{\pi(2i-1)!}\left( \prod_{j=1}^{i-1} (j^2+ \frac{(\log|x|)^2}{\pi^2}) \right)\frac{\log|x|}{x^2-1}.
\end{align*}
Furthermore, one might use the ``universal prior" due to \cite{rissanen1983universal} over the number of terms $k$ in the product density. The ``universal prior" is defined with the mass function: 
$$
Q(i) = 2^{-L^0(i)}, \mbox{ for } i=1, 2, \ldots; \quad L^0(i) = \log^*(i) + \log c,
$$
where, $\log^*(x) = \log x + \log \log x + \ldots$, where the sum involves only non-negative terms and $c = \sum 2^{-\log^*i} \approx 2.865064$. \par
The family of Cauchy-product densities can be used in conjunction with Rissanen's universal prior described above to define an adaptive shrinkage estimator such as the \textit{Polyshrink} estimator due to \cite{foster2005polyshrink}, where the amount of shrinkage varies adaptively with the estimation task. For an $n$-dimensional parameter $\theta$, the Polyshrink estimator uses a collection of discrete mixture models $\mathcal{G}_{p} = \left\{ g_{\epsilon_k}(y) = (1-\epsilon_k)\phi(y)+\epsilon_k \psi(z), \; \epsilon_k = 2^{k-(K+1)} \right\}$, for $k = 1, \ldots, K$ with $K=1+ \lfloor{\log_2(p)}\rfloor$ and $\phi(\cdot)$ and $\psi(\cdot)$ denote the standard normal density and Cauchy density with scale $\sqrt{2}$ respectively. We conjecture the advantages of the one-group model over the two-groups model would naturally carry over to this case if we use a collection of one-group priors defined by Cauchy products of different orders to achieve different amounts of shrinkage. 
The possibility of such extensions was first discussed in \cite{polson2012half}, and it would be interesting to settle this issue theoretically.

\appendix
\section{Proofs of theorems}
\renewcommand{\theequation}{A.\arabic{equation}}
\renewcommand{\thesubsection}{A.\arabic{subsection}}

\subsection{Proof of Theorem \ref{th:1}} \label{app:th1}
Let $p_{HS+}( \theta )$ be the marginal density for the horseshoe+ prior. From Equations (\ref{eq:hs+-hier}-\ref{eq:lambdai}) with $\tau=1$, we have,
$$
p_{HS+}( \theta ) =  \frac{4}{\pi^2} \frac{1}{\sqrt{2 \pi}} \int_0^\infty \frac{1}{\lambda} e^{- \frac{1}{2 \lambda^2} \theta^2} \frac{ \log ( \lambda ) }{ \lambda^2 -1 } d \lambda.
$$
First note that,  applying the transformation $\zeta = 1/\lambda^2$ we get 
\begin{equation*}
p_{HS+}(\theta)  = \frac{1}{\pi^2 \sqrt{2\pi}}\int_0^\infty e^{- \frac{\zeta}{2} \theta^2} \frac{ \log ( \zeta) }{ \zeta -1 } d \zeta.
\end{equation*}
For an upper bound one can use the fact that 
$$
\frac{\log \zeta}{\zeta-1} \leq \frac{1}{\sqrt{\zeta}}\quad \text{for} \quad \zeta>0. 
$$
This gives
\begin{align*}
\int_0^\infty e^{- \frac{\zeta}{2} \theta^2} \frac{ \log ( \zeta ) }{ \zeta -1 } d \zeta & \leq \int_0^\infty  \zeta^{-1/2} e^{- \frac{\zeta}{2} \theta^2} d \zeta\\
&= \frac{\Gamma(1/2)}{(\theta^2/2)^{1/2}} = \frac{\sqrt{2 \pi}}{|\theta|}.
\end{align*}
For a lower bound one can use the fact that 
$$
\frac{\log \zeta}{\zeta-1} \geq \frac{2}{1 + \zeta}\quad \text{for} \quad \zeta>0. 
$$
This gives
\begin{align*}
\int_0^\infty e^{- \frac{\zeta}{2} \theta^2} \frac{ \log ( \zeta ) }{ \zeta -1 } d \tau &\geq \int_0^\infty e^{- \frac{\zeta}{2} \theta^2} \frac{2}{ \zeta +1 } d \zeta\\
& = 2 e^{\theta^2/2} E_1(\theta^2/2) \\
&> \log(1 + \frac{4}{\theta^2}),
\end{align*}
where $E_1(\cdot)$ is the exponential integral. Thus, combining the lower and upper bounds,
$$
\frac{1}{\pi^2 \sqrt{2\pi}} \log (1 + \frac{4}{\theta^2}) < p_{HS+} (\theta) \leq  \frac{1}{\pi^2 |\theta|}.
$$
This completes the proof of Part (1). Part (2) then follows from the fact that the lower bound goes to $\infty$ as $\theta$ goes to zero. In comparison, horseshoe has bounds
$$\frac{K}{2}\log (1 + \frac{4}{\theta^2}) < p_{HS}(\theta) < K \log (1 + \frac{2}{\theta^2}),$$
with $K=1/(2\pi^3)^{1/2}$. Both the upper and lower bounds are sharper compared to horseshoe+. However, the bounds for horseshoe+ can be sharpened using better approximations for $\log(\zeta)/(\zeta-1)$, for which an infinite product representation is given by
$$
\frac{\log \zeta}{\zeta-1} = \frac{2}{1 + \sqrt{\zeta}} \times \frac{2}{1 + \sqrt{\sqrt{\zeta}}} \times \ldots = \prod_{i=1}^{\infty} \frac{2}{1 + \zeta^{1/2^i}}.
$$

\subsection{Proof of Theorem \ref{th:3}} \label{app:th3}

We write the posterior density of $\kappa_i$ given $y_i$ in $\tau$ scale as follows: 
$$
p(\kappa_i | y_i, \tau) \propto (1-\kappa_i)^{-\half} \exp \left \{ - \kappa_i \frac{ y_i^2}{2} \right \} \frac{|\log \left \{ \kappa_i \tau^2 / ( 1 - \kappa_i ) \right \} |}{|(\kappa_i (\tau^2 +1 ) -1)|}.
$$
Now, we use the inequality $ 1- 1/x < \log (x) < x-1$ for $x > 0$ to derive the following bounds for the Jacobian term:
\begin{align*}
1-\frac{1-\kappa_i}{\kappa_i \tau^2} & < |\log \left \{ \kappa_i \tau^2 / ( 1 - \kappa_i ) \right \} | < \frac{\kappa_i \tau^2}{( 1 - \kappa_i )} - 1, \\
\frac{1}{\kappa_i \tau^2} & < \frac{|\log \left \{ \kappa_i \tau^2 / ( 1 - \kappa_i ) \right \} |}{|(\kappa_i (\tau^2 +1 ) -1)|} < \frac{1}{1-\kappa_i}.
\end{align*}
First note that $\P(\kappa_i < \epsilon | y_i, \tau) \leq  \P(\kappa_i < \epsilon | y_i, \tau)/\P(\kappa_i > \epsilon | y_i, \tau)$. Moreover, we can bound the ratio $\P(\kappa_i < \epsilon | y_i, \tau)/\P(\kappa_i > \epsilon | y_i, \tau)$ as follows: 
\begin{align*} 
\frac{\P(\kappa_i < \epsilon | y_i, \tau)}{\P(\kappa_i > \epsilon | y_i, \tau)} & = \frac{\int_{0}^{\epsilon} (1-\kappa_i)^{-\half} \exp \left \{ - \kappa_i \frac{ y_i^2}{2} \right \} \frac{|\log \left \{ \kappa_i \tau^2 / ( 1 - \kappa_i ) \right \} |}{|(\kappa_i (\tau^2 +1 ) -1)|} d\kappa_i }{\int_{\epsilon}^{1} (1-\kappa_i)^{-\half} \exp \left \{ - \kappa_i \frac{ y_i^2}{2} \right \} \frac{|\log \left \{ \kappa_i \tau^2 / ( 1 - \kappa_i ) \right \} |}{|(\kappa_i (\tau^2 +1 ) -1)|} d\kappa_i} \\
& \leq \frac{\int_{0}^{\epsilon} (1-\kappa_i)^{-\half} \exp \left \{ - \kappa_i \frac{ y_i^2}{2} \right \} (1-\kappa_i)^{-1} d\kappa_i }{\int_{\epsilon}^{1} (1-\kappa_i)^{-\half} \exp \left \{ - \kappa_i \frac{ y_i^2}{2} \right \} (\kappa_i \tau^2)^{-1} d\kappa_i} \\
& \leq \frac{\int_{0}^{\epsilon} (1-\kappa_i)^{-3/2} d\kappa_i }{\exp \left \{ -\frac{ y_i^2}{2} \right \}\frac{1}{\tau^2} \int_{\epsilon}^{1} (1-\kappa_i)^{-\half} \frac{1}{\kappa_i} d\kappa_i} \\
& \leq \mathrm{e}^\frac{ y_i^2}{2}  \tau^2 \frac{(1-\epsilon)^{-3/2} \epsilon}{1/\sqrt{1-\epsilon}} = \mathrm{e}^\frac{ y_i^2}{2} \tau^2 \epsilon (1-\epsilon)^{-2}.
\end{align*}
The final step follows from the penultimate step by bounding the integral by the extreme values of the integrand multiplied by the length of the interval of integration. 
\qedhere


\subsection{Proof of Theorem \ref{th:typeI}} \label{app:typeI}
We need the following inequality proved in Appendix~\ref{app:th3}: 
\begin{equation*}
\frac{1}{\kappa_i \tau^2}  < \frac{|\log \left \{ \kappa_i \tau^2 / ( 1 - \kappa_i ) \right \} |}{|(\kappa_i (\tau^2 +1 ) -1)|} < \frac{1}{1-\kappa_i}. \label{t1-1}
\end{equation*}
To use the probability concentration inequality in Theorem \ref{th:3} to prove the bound on type-I error rate, we need to establish the following fact:
\begin{equation*}
\E(\kappa_i | y_i, \tau) = \P(\kappa_i > \half | y_i, \tau)(1+o(1)) \label{t1-2}.
\end{equation*}
We will prove this in two steps. We first show that the posterior mean can be well approximated by evaluating the integral from $\half$ to $1$, i.e. 
$$
\E(\kappa_i | y_i,\tau) = \int_{\half}^{1} \kappa_i p(\kappa_i | y_i, \tau) d \kappa_i (1+o(1)),
$$ 
as $\tau \rightarrow 0$. In the next step, we prove that 
$$
\int_{\half}^{1} \kappa_i p(\kappa_i | y_i, \tau) d \kappa_i  = \int_{\half}^{1} p(\kappa_i | y_i, \tau) d \kappa_i (1+o(1)),
$$ 
as $\tau \rightarrow 0$. First note that: 
\begin{align*}
\frac{\int_{0}^{\half} \kappa_i p(\kappa_i | y_i, \tau) d \kappa_i }{\int_{\half}^{1} \kappa_i p(\kappa_i | y_i, \tau) d \kappa_i} & \leq 
\frac{ e^{\frac{y_i^2}{2}} \int_{0}^{\half} \kappa_i (1-\kappa_i)^{-\half} (1-\kappa_i)^{-1} d \kappa_i}{\int_{\half}^{1} \kappa_i (1-\kappa_i)^{-\half} \frac{1}{\kappa_i\tau^2} d \kappa_i}. \\
& = e^{\frac{y_i^2}{2}} \tau^2 \frac{ \int_{0}^{\half} \kappa_i (1-\kappa_i)^{-3/2} d \kappa_i}{\int_{\half}^{1} (1-\kappa_i)^{-\half} d \kappa_i} \\
& \leq (3-2\sqrt{2}) e^{\frac{y_i^2}{2}} \tau^2. 
\end{align*}
Thus,
$$
\E(\kappa_i | y_i,\tau) = \int_{0}^{\half} \kappa _i p(\kappa_i | y_i, \tau) d \kappa_i +   \int_{\half}^{1} \kappa_i p(\kappa_i | y_i, \tau) d \kappa_i =\int_{\half}^{1} \kappa_i p(\kappa_i | y_i, \tau) d \kappa_i (1+o(1))
$$
as $\tau \to 0$. Next, note that
\begin{align*}
\frac{\int_{\half}^{1} (1-\kappa_i) p(\kappa_i | y_i, \tau) d \kappa_i}{\int_{\half}^{1} p(\kappa_i | y_i, \tau) d \kappa_i} & \leq \frac{ e^{\frac{y_i^2}{2}} \int_{\half}^{1} (1-\kappa_i)^{-\half} d \kappa_i}{\int_{\half}^{1} (1-\kappa_i)^{-\half} \frac{1}{\kappa_i\tau^2} d \kappa_i} \\
& \leq e^{\frac{y_i^2}{2}} \tau^2. 
\end{align*}
Thus, we have $\int_{\half}^{1} \kappa_i p(\kappa_i |y_i, \tau) d \kappa_i = \int_{\half}^{1} p(\kappa_i |y_i, \tau) d \kappa_i (1 -o(1))$. 
The asymptotic expression for the type-I error rate is then calculated as 
\begin{align*}
\P_{H_0} \left( \E(\kappa_i | y_i, \tau) < \frac{1}{2}\right ) &= \P_{H_0}\left( \int_{\half}^{1} p(\kappa_i | y_i, \tau) d\kappa_i< \half \right) (1+o(1))\\
&= \P_{H_0}\left( \int_{0}^{\half} p(\kappa_i | y_i, \tau) d\kappa_i > \half \right)(1+o(1)). 
\end{align*}
Applying Theorem \ref{th:3} for $\epsilon = \half$, we have $\int_{0}^{\half} p(\kappa_i | y_i, \tau) d\kappa_i \leq (2\mathrm{e}^\frac{ Y_i^2}{2} \tau^2)(1+o(1))$. 
\begin{align*}
\P_{H_0} \left( \E(\kappa_i | y_i, \tau) < \half \right ) & \leq \P_{Y_i \sim \Nor(0,1)}( 2\mathrm{e}^\frac{ Y_i^2}{2} \tau^2 > \half) (1+o(1)) \\
&= \P(|Y_i| > 2 \sqrt{\log(1/2\tau)}) (1+o(1)) \\
& \approx \frac{\phi(2 \sqrt{\log(1/2\tau)})}{2 \sqrt{\log(1/2\tau)}}(1+o(1)) \\
& = \sqrt{\frac{2}{\pi}}\frac{\tau^2}{\sqrt{\log(1/2\tau)}}(1+o(1)).
\end{align*}


\subsection{Proof of Theorem \ref{th:5}} \label{app:th5}

First note the following identity:
\begin{align*}
\frac{1-\kappa_i}{\kappa_i \tau^2} & = \left \{ \frac{(1+\tau^2)(1-\kappa_i)}{\tau^2} \times \frac{1}{(1+\tau^2)\kappa_i} \right \} \\ 
& =  \left \{ (1 + \frac{(1+\tau^2)(1-\kappa_i)-\tau^2}{\tau^2}) \times (1+ \frac{1/(1+\tau^2)-\kappa_i}{\kappa_i}) \right \} \\
& = \left \{ (1 + \frac{(1-\kappa_i(1+\tau^2))}{\tau^2}) \times (1+ \frac{(1-\kappa_i(1+\tau^2))}{\kappa_i(1+\tau^2)}) \right \} .
\end{align*}
Since both $\frac{(1+\tau^2)(1-\kappa_i)-\tau^2}{\tau^2} \geq -1$ and $\frac{1/(1+\tau^2)-\kappa_i}{\kappa_i} \geq -1$, we get the following upper bound to the Jacobian term of the horseshoe+ posterior density: 
\begin{align}
\frac{|\log (\frac{1-\kappa_i}{\kappa_i \tau^2})|}{|1-\kappa_i(1+\tau^2)|} & = \frac{|\log \left(1 + \frac{(1-\kappa_i(1+\tau^2))}{\tau^2} \right) + \log \left(1+ \frac{(1-\kappa_i(1+\tau^2))}{\kappa_i(1+\tau^2)} \right) |}{|1-\kappa_i(1+\tau^2)|} \nonumber\\
& \leq \frac{1}{\tau^2} + \frac{1}{\kappa_i(1+\tau^2)}. \label{upper_jac}
\end{align}
Also, for $\kappa_i < 1/(1+\tau^2)$, we would have $\log(\frac{1-\kappa_i}{\kappa_i \tau^2}) \geq \frac{1-\kappa_i(1+\tau^2)}{1-\kappa_i}$. 
Now we use the upper bound in Equation \eqref{upper_jac} and the above lower bound to derive an upper bound to the tail probability of the posterior distribution of $\kappa_i$ given $\tau$ and $y_i$. We will use the fact that
$$
\P(\kappa_i > \eta | y_i, \tau) \leq \P(\kappa_i > \eta | y_i, \tau)/\P(\kappa_i < \eta | y_i, \tau) \leq \P(\kappa_i > \eta | y_i, \tau)/\P(\kappa_i < \eta \delta | y_i, \tau),
$$
for some fraction $\delta \in (0,1)$,  such that $\eta\delta\leq1/(1+\tau^2)$. Thus,
\begin{align*}
\P(\kappa_i > \eta | \tau, y_i) & \leq \frac{\int_{\eta}^{1} \frac{1}{\sqrt{1-\kappa_i}} \exp \left \{ - \kappa_i \frac{ y_i^2}{2} \right \}
 \frac{|\log (\frac{1-\kappa_i}{\kappa_i \tau^2})|}{|1-\kappa_i(1+\tau^2)|}  d\kappa_i} {\int_{0}^{\eta\delta} \frac{1}{\sqrt{1-\kappa_i}} \exp \left \{ - \kappa_i \frac{ y_i^2}{2} \right \} \frac{|\log (\frac{1-\kappa_i}{\kappa_i \tau^2})|}{|1-\kappa_i(1+\tau^2)|}  d\kappa_i} \\ 
& \leq e^{-\eta(1-\delta)\frac{y_i^2}{2}} \frac{\int_{\eta}^{1} \frac{1}{\sqrt{1-\kappa_i}} \left(\frac{1}{\tau^2} + \frac{1}{\kappa_i(1+\tau^2)} \right) d\kappa_i} {\int_{0}^{\eta\delta} \frac{1}{\sqrt{1-\kappa_i}} \frac{1}{1-\kappa_i} d\kappa_i} \\
& \leq e^{-\eta(1-\delta)\frac{y_i^2}{2}} \frac{ \left\{ \frac{\sqrt{1-\eta}}{\tau^2} + \frac{\mathrm{arctanh} (\sqrt{1-\eta})}{(1+\tau^2)} \right\}} { \{ 1/\sqrt{1-\eta\delta} -1 \} } \\
& \leq e^{-\eta(1-\delta)\frac{y_i^2}{2}} \frac{1}{\tau^2}\frac{ \left\{\sqrt{1-\eta} + \mathrm{arctanh} (\sqrt{1-\eta}) \right\}} { \{ 1/\sqrt{1-\eta\delta} -1 \} }\\
& = e^{-\eta(1-\delta)\frac{y_i^2}{2}} \frac{1}{\tau^2} C(\eta, \delta),
\end{align*}
where, 
$$C(\eta, \delta) = \frac{ \left\{\sqrt{1-\eta} + \mathrm{arctanh} (\sqrt{1-\eta}) \right\}} { \{ 1/\sqrt{1-\eta\delta} -1 \} },$$
is independent of $y_i$.
\qedhere

\subsection{Proof of Theorem \ref{th:typeII}} \label{app:typeII}
The probability of type-II error for any global-local shrinkage rule is given by: 
\begin{equation*}
t_{2} = \P_{Y_i \sim f_{H_A}(y)}\left( \mathbb{E}(\kappa_i | y_i, \tau) > \half \right). 
\end{equation*}
We note that, from Assumption \ref{assnA} and the aforementioned choice of $\tau$, we have 
$
\log(1/\tau^2)/\psi^2 \rightarrow C \in (0,\infty),
$
where $C$ is the threshold appearing in the expression for the risk of Bayes oracle (by Equation \eqref{optrisk}). First note that, for any $\eta \in (0,1)$, 
\begin{align*}
\kappa_i & < I( \kappa_i \in (\eta,1)) + \eta, \\
\mathbb{E} (\kappa_i | y_i,\tau) & < \P(\kappa_i > \eta | y_i, \tau) + \eta.
\end{align*} 
Then the type-II error rate can be written as: 
\begin{align*}
t_2 & = \P_{Y_i \sim \Nor(0,1+\psi^2)} \left( \mathbb{E} (\kappa_i | y_i, \tau) > \half \right) \\
& \leq \P_{Y_i \sim \Nor(0,1+\psi^2)} \left( \P_{[\kappa_i |y_i ,\tau]}(\kappa_i > \eta | y_i, \tau) > \half-\eta \right) 
\end{align*}
We use the upper bound $\P(\kappa_i > \eta | y_i, \tau) \leq \exp(-\eta(1-\delta)Y_i^2/2)\frac{1}{\tau^2} C(\eta,\delta)$ from the concentration inequality in Theorem \ref{th:5} to derive the following: 

\begin{align*}
t_2 & \leq \P_{Y_i \sim \Nor(0,1+\psi^2)} \left( \exp(-\eta(1-\delta)Y_i^2/2)\frac{1}{\tau^2} C(\eta,\delta) > \half-\eta \right) \\
& \leq \P_{Y_i \sim \Nor(0,1+\psi^2)} \left( Y_i^2 < \frac{2}{\eta(1-\delta)} \left\{ \log(\frac{C(\eta,\delta)}{\half-\eta}) + \log(\frac{1}{\tau^2}) \right\} \right) \\
& \leq \P_{Z_i  \sim \Nor(0,1)} \left( |Z_i| < \sqrt{\frac{2}{\eta(1-\delta)}} \sqrt{\frac{\log(1/\tau^2)}{\psi^2}}(1+o(1)) \right) \mbox{ as } n \rightarrow \infty.
\end{align*}
Under the assumption $\tau = O(\mu)$, we have $\frac{\log(1/\tau^2)}{\psi^2} \rightarrow C$ as $n \rightarrow \infty$, where $C$ is the constant appearing in the Bayes risk for the oracle (by Equation \eqref{optrisk}). Thus, we have:

\begin{align*}
t_2 \leq  \left( 2\Phi(\sqrt{\frac{2}{\eta(1-\delta)}} \sqrt{C}) - 1 \right) (1+ o(1)) \mbox{ as } n \rightarrow \infty.
\end{align*}


\qedhere
\subsection{Proof of Theorem \ref{th:int}}\label{app:meijer}
From Appendix~\ref{app:th1}, the marginal prior, $p_{HS+}(\theta)$, for horseshoe+ is given by the convolution
\begin{align}
  p_{HS+}(\theta) &= \frac{1}{\sqrt{2}\pi^{5/2}}\int_0^\infty e^{- \frac{\zeta}{2} \theta^2}
  \frac{ \log(\zeta) }{ \zeta -1 } d\zeta. 
  \label{eq:hs-plus-marginal}
\end{align}
We can use Meijer G-function convolution identities \citep{mathai2009h}
to provide a special function representation of such distributions.
Let $\MeijerG*{m}{n}{p}{q}{\ba_p}{\bb_q}{x}$ denote the Meijer-G function,
This class of functions satisfies a very useful convolution identity with many
applications in Bayesian computation. We have the integral identity
\begin{align*}
  \int_0^\infty  
    \MeijerG*{m}{n}{p}{q}{\ba_p}{\bb_q}{\eta x} 
    \MeijerG*{\mu}{\nu}{\sigma}{\tau}{\bc_\sigma}{\bd_\tau}{\omega x} dx
    &= \frac{1}{\eta} \MeijerG*{n+\mu}{m+\nu}{q+\sigma}{p+\tau}{-b_1,\dots,-b_m,\bc_\sigma,-b_{m+1},\dots,-b_q}{-a_1,\dots,-a_n,\bd_\tau, -a_{n+1}, \dots, -a_p}{\frac{\omega}{\eta}} \\
    &= \frac{1}{\omega} \MeijerG*{m+\nu}{n+\mu}{p+\tau}{q+\sigma}{a_1,\dots,a_n,-\bd_\tau,a_{n+1},\dots,a_p}{b_1,\dots,b_m,-\bc_\sigma, b_{m+1}, \dots, b_q}{\frac{\eta}{\omega}} .
\end{align*}
Furthermore, when one of the G-functions is an exponential function, we obtain a 
general form for the Laplace transform of a G-function:
\begin{align*}
  \int_0^\infty  
    e^{-\omega x}
    \MeijerG*{m}{n}{p}{q}{\ba_p}{\bb_q}{\eta x} dx
    &= \omega^{-1} 
    \MeijerG*{m}{n+1}{p+1}{q}{0, \ba_p}{\bb_q}{\eta \omega^{-1}},
\end{align*}
for $\operatorname{Re}(s) > 0$.
This can be used to calculate the implied prior via the following identity:
\begin{align*}
  \frac{\log(x)}{x-1} &= \MeijerG*{2}{2}{2}{2}{0}{0}{x}, 
\end{align*}
for appropriate $x$. Equation~\eqref{eq:hs-plus-marginal} is then
\begin{align*}
  p_{HS+}(\theta) =&\frac{1}{\sqrt{2}\pi^{5/2}}
  \MeijerG*{2}{3}{3}{2}{1,1,1}{1,1}{2 \theta^{-2}}. 
\end{align*}
This identity allows us to develop the asymptotic behavior around 
$\theta \to 0$ and $\infty$ through power-logarithmic series expansions
for G-functions \citep{kilbas1900h, mathai2009h}.  
The dominant terms in a power series expansion around $\infty$ are
\begin{align}
  p^{\infty}_{HS+}(\theta) &= \frac{1}{\sqrt{2} \pi^{5/2}} 
    \frac{2 \log(\theta)+\gamma -\log(2)}{\theta^2} \;, 
    \label{eq:hs-plus-asymp-inf}
\end{align}
where $\gamma$ is the Euler-Mascheroni constant, and similarly around zero
\begin{multline}
  p^{0}_{HS+}(\theta) = 
  \frac{1}{24 \sqrt{2} \pi^{5/2}} \left(
    24 \log^2\left(\frac{1}{\theta}\right)+24 \log(2) \log\left(\frac{1}{\theta}\right) \right. \\
    \left. -24 \gamma \log\left(\frac{1}{\theta}\right)+6 \gamma^2
    +5 \pi^2+6 \log^2(2)-12 \gamma \log(2) \right).
    \label{eq:hs-plus-asymp-origin}
\end{multline}
The corresponding results for the horseshoe prior are 
obtained via the substitution $u=1/\lambda^2$ (as in \cite{carvalho2010horseshoe})
and the identity
\begin{align*}
  p_{HS}(\theta) &\propto \int_0^\infty \frac{1}{1+u} e^{- \frac{\theta^2 u}{2}} du \\
                   &= \MeijerG*{1}{2}{2}{1}{1,1}{1}{2 \theta^{-2}} ,
\end{align*}
from which we obtain respectively around $\infty$ and $0$ that
\begin{align}
  p^{\infty}_{HS}(\theta) &= \frac{\sqrt{2}}{\pi^{3/2}} \frac{1}{\theta^2},
    \label{eq:hs-asymp-inf}
\end{align}
and
\begin{align}
  p^{0}_{HS}(\theta) &= \frac{1}{\sqrt{2} \pi^{3/2}}\left(
      2 \log\left(\frac{1}{\theta}\right)-\gamma +\log(2) \right). \;
    \label{eq:hs-asymp-origin}
\end{align}
 The behavior of the Bayes risk depends crucially on the amount of
mass the prior places around the origin.
Using Equation \eqref{eq:hs-plus-asymp-origin} and integrating (on one side) around zero,
\begin{multline*}
  \int_0^{\frac{1}{\sqrt{n}}} p^0_{HS+}(\theta) d\theta =
  \frac{1}{24 \sqrt{2} \pi^{5/2} \sqrt{n}} 
  \left( 6 \log(n) (\log(n)-2 \gamma +4+\log(4))+6 \gamma^2+5 \pi^2
    \right.\\
    +6 \left. \left(8+\log^2(2)+\log(16)\right)-12 \gamma (2+\log(2)) \right).
\end{multline*}
Collecting the higher order terms in $n$ we get
\begin{align*}
 \frac{1}{\sqrt{2} \pi^{5/2} \sqrt{n}} \left(
    \frac{\log^2(n)}{4}+
    \left(1 -\frac{\gamma}{2}
    +\frac{\log(4)}{4}
    \right) \log(n) + O(1) 
  \right).  \; 
\end{align*}
Similarly, using Equation (\ref{eq:hs-asymp-origin}) we have
$$
\int_{0}^{\frac{1}{\sqrt{n}}} p^0_{HS} (\theta) d \theta = \frac{1}{\sqrt{2} \pi^{3/2} \sqrt{n}} \left( \frac{\log(n)}{2}  + O(1) \right),
$$
completing the proof of Theorem~\ref{th:int}.

\subsection{Proof of Theorem \ref{th:mse}} \label{app:mse}
Using Equation \eqref{hseq1} in Equation \eqref{tweedie1} for large values of $|y_i|$ yields,
\begin{eqnarray*}
\E_{HS+} (\theta_i|y_i)  &=& y_i + \frac{d}{dy_i} [\log m_{HS+}(y_i)]\\
&=& \E_{HS} (\theta_i|y_i) + \frac{1}{y_i\log |y_i|} - O\left(\frac{1}{y_i^2}\right).
\end{eqnarray*}
Thus,
$$
\mathrm{Bias}_{HS+} (\theta_i|y_i) = \mathrm{Bias}_{HS} (\theta_i|y_i) + \frac{1}{y_i\log |y_i|} - O\left(\frac{1}{y_i^2}\right).
$$
Continuing the calculation for variances we have (by Equation \eqref{hseq2} and \eqref{tweedie2}),
\begin{eqnarray*}
\V_{HS+} (\theta_i|y_i)  &=& 1 + \frac{d^2}{dy_i^2} [\log m_{HS+}(y_i)]\\
&=& \V_{HS}(\theta_i| y_i)  - \frac{1+ \log |y_i|}{(y_i\log |y_i|)^2} + O\left(\frac{1}{y_i^3}\right).
\end{eqnarray*}
Thus, comparing the posterior MSE of the two estimators for large values of $y$ we have
\begin{eqnarray*}
\MSE_{HS+} (\theta_i|y_i) &=& \mathrm{Bias}^2_{HS+} (\theta_i|y_i) + \V_{HS+} (\theta_i|y_i) \\
&=& (\mathrm{Bias}_{HS} (\theta_i|y_i) + \frac{1}{y_i\log |y_i|} + O(1/y_i^2))^2 + \V_{HS}(\theta_i|y_i)  - \frac{1+ \log |y_i|}{(y_i\log |y_i|)^2} + O\left(\frac{1}{y_i^3}\right) \\
&=& \MSE_{HS} (\theta_i|y_i) + 2\frac{\mathrm{Bias}_{HS} (\theta_i|y_i)}{y_i\log |y_i|} - \frac{1}{y_i^2\log |y_i|} + O \left(\frac{1}{y_i^3\log |y_i|}\right) +O\left(\frac{1}{y_i^3}\right)\\
&=& \MSE_{HS} (\theta_i|y_i) - \frac{1}{y_i^2\log |y_i|} + O \left(\frac{1}{y_i^3}\right) ,
\end{eqnarray*}
where the last line follows from the fact that $\mathrm{Bias}_{HS} (\theta_i|y_i) = O(1/y_i^2)$ (from Theorem 3 in \cite{carvalho2010horseshoe}). \qedhere


\bibliographystyle{biom}
\bibliography{horseshoe-plus}

\end{document}